\documentclass[a4paper, twoside, 12pt,reqno]{amsart}

\usepackage[left=3cm,right=4cm,top=3cm,bottom=3.5cm,footskip=1cm]{geometry}

\usepackage{amssymb}
\usepackage{amsmath}
\usepackage{amsfonts}
\usepackage{amsthm}
\usepackage{enumerate}

\usepackage[latin1]{inputenc}
\usepackage[T1]{fontenc}

\usepackage[english]{babel} 

\usepackage{graphicx}

\usepackage{ifpdf}
\ifpdf
\fi

\usepackage{bbm}

\usepackage[mathscr]{euscript}

\renewcommand{\mathcal}{\mathscr}

\usepackage{ifthen}

\usepackage[all,line,ps]{xy}
\SelectTips{cm}{12}
\CompileMatrices

\usepackage{url}

\newcommand{\mathset}[1]{\mathbbm{#1}}

\newcommand{\setN}{\mathset{N}}
\newcommand{\setZ}{\mathset{Z}}

\newcommand{\setR}{\mathset{R}}
\newcommand{\setC}{\mathset{C}}

\newcommand{\G}[1]{\mathit{MT}(#1)}
\newcommand{\RP}{\mathset{R}P}
\newcommand{\CP}{\mathset{C}P}

\newcommand{\CAT}{\mathbf{CAT}}

\DeclareMathOperator{\pr}{pr}
\DeclareMathOperator*{\colim}{colim}
\DeclareMathOperator*{\hocolim}{hocolim}
\DeclareMathOperator*{\hofib}{hofib}

\DeclareMathOperator{\Emb}{Emb}
\DeclareMathOperator{\Diff}{Diff}
\DeclareMathOperator{\Spin}{Spin}

\DeclareMathOperator{\Map}{Map}

\DeclareMathOperator{\id}{id}
\DeclareMathOperator{\Th}{Th}
\DeclareMathOperator{\Int}{Int}
\DeclareMathOperator{\GL}{GL}
\DeclareMathOperator{\Path}{Path}

\DeclareMathOperator{\ob}{ob}
\DeclareMathOperator{\mor}{mor}

\DeclareMathOperator{\supp}{supp}
\DeclareMathOperator{\proj}{proj}
\DeclareMathOperator{\sgn}{sgn}
\DeclareMathOperator{\Bun}{Bun}

\newcommand{\cdp}{C_d^\pitchfork}

\newcommand{\ddp}{D_d^\pitchfork}
\newcommand{\ddpb}{D_{d,\partial}^\pitchfork}

\newcommand{\udnb}{U_{d, n}^\bot}
\newcommand{\point}{\mathrm{pt}}
\newcommand{\ShX}{\mathrm{Sh}(\mathcal{X})} 
\newcommand{\isom}{\cong}

\renewcommand{\phi}{\varphi}
\renewcommand{\epsilon}{\varepsilon}

\newcommand{\inv}{^{-1}}

\newcommand{\abs}[2][]{#1\lvert #2#1\rvert}

\numberwithin{equation}{section}

\theoremstyle{plain}
\newtheorem{lemma}{Lemma}[section]
\newtheorem{proposition}[lemma]{Proposition}
\newtheorem{theorem}[lemma]{Theorem}
\newtheorem{criteria}[lemma]{Criteria}
\newtheorem{corollary}[lemma]{Corollary}
\newtheorem{construction}[lemma]{Construction}

\newtheorem*{mthm}{Main Theorem}
\newtheorem*{mthm2}{Main Theorem (with tangential structures)}

\theoremstyle{definition}
\newtheorem{definition}[lemma]{Definition}
\newtheorem{remark}[lemma]{Remark}

\newenvironment{iblist}{
  \begin{list}{(\roman{enumi})\hfil }{
      \labelsep=10pt
      \itemindent=10pt
      \leftmargin=0pt
      \labelwidth=0pt
      \usecounter{enumi}
    }}{\end{list}}

\author{S. Galatius, I. Madsen, U.  Tillmann and M. Weiss}

\title{The homotopy type of the cobordism category}

\date{\today}
\thanks{S.\ Galatius is partially supported by NSF grant DMS-0505740 and the Clay Institute}
\thanks{M.\ Weiss is partially supported by the Royal Society and
  the Mittag-Leffler Institute}
\address{Stanford University, Stanford, California}
\email{galatius@math.stanford.edu}
\address{University of Copenhagen, Copenhagen, Denmark}
\email{imadsen@math.ku.dk}
\address{University of Oxford, Oxford, United Kingdom}
\email{tillmann@maths.ox.ac.uk}
\address{University of Aberdeen, Aberdeen, United Kingdom}
\email{m.weiss@maths.abdn.ac.uk}

\begin{document}

\begin{abstract}
The embedded cobordism category under study in this paper generalizes
the category of conformal surfaces, introduced by G.\ Segal in
\cite{MR2079383} in order to formalize the concept of field theories.
Our main result identifies the homotopy type of the classifying space
of the embedded $d$-dimensional cobordism category for all $d$.  For
$d=2$, our results lead to a new proof of the generalized Mumford
conjecture, somewhat different in spirit from the original one,
presented in \cite{arXiv1}.
\end{abstract}

\maketitle

\tableofcontents

\section{Introduction and results}
\label{cha:introduction-results}

The conformal surface category $\mathscr{S}$ is defined as
follows.  For each non-negative integer $m$ there is one object $C_m$
of $\mathscr{S}$, namely the 1-manifold $S^1 \times \{1,2, \dots,
m\}$.  A morphism from $C_m$ to $C_n$ is an isomorphism class of a
Riemann surface $\Sigma$ with boundary $\partial \Sigma$ together with
an orientation-preserving diffeomorphism $\partial \Sigma \to C_n
\amalg - C_m$.  The composition is by sewing surfaces together.

Given a differentiable subsurface $F \subseteq [a_0, a_1] \times
\setR^{n+1}$ with $\partial F = F \cap \{a_0, a_1\} \times
\setR^{n+1}$, each tangent space $T_p F$ inherits an inner product
from the surrounding euclidean space and hence a conformal structure.
If $F$ is oriented, this induces a complex structure on $F$.
Associating a complex structure to an embedded surface in this way is,
suitably interpreted, a homotopy equivalence (namely the space of
complex structures and the space of embeddings in $\setR^\infty$ are
both contractible.  See Remark~\ref{remark:Segal} for a further
discussion).  The category $\mathcal{C}_2$ of embedded oriented
surfaces can thus be viewed as a substitute for the conformal surface
category.

The embedded surface category has an obvious generalization to higher
dimensions.  For any $d \geq 0$, we have a category $\mathcal{C}_d$
whose morphisms are $d$-dimensional submanifolds $W \subseteq [a_0,
a_1] \times \setR^{n+d-1}$ that intersect the walls $\{a_0,a_1\}
\times \setR^{n+d-1}$ transversely in $\partial W$.  The codimension
$n$ is arbitrarily large, and not part of the structure.  Viewing $W$
as a morphism from the incoming boundary $\partial_\mathrm{in} W =
\{a_0\} \times \setR^{n+d-1} \cap W$ to the outgoing boundary
$\partial_\mathrm{out} W = \{a_1\} \times \setR^{n+d-1} \cap W$, and
using union as composition, we get the embedded cobordism category
$\mathcal{C}_d$.

It is a topological category in the sense that the total set of
objects and the total set of morphisms have topologies such that the
structure maps (source, target, identity and composition) are
continuous.  In fact, there are homotopy equivalences
\begin{align*}
  \ob\mathcal{C}_d \simeq \coprod_M B \Diff(M), \quad
  \mor\mathcal{C}_d \simeq \coprod_W B \Diff(W;
  \{\partial_\mathrm{in}\}, \{\partial_\mathrm{out}\})
\end{align*}
where $M$ varies over closed $(d-1)$-dimensional manifolds and $W$
over $d$-dimensional cobordisms, one in each diffeomorphism class.
Here $\Diff(M)$ denotes the topological group of diffeomorphisms of
$M$ and $\Diff(W, \{\partial_\mathrm{in}\},
\{\partial_\mathrm{out}\})$ denotes the group of diffeomorphisms of
$W$ that restrict to diffeomorphisms of the incoming and outgoing
boundaries.  Source and target maps are induced by restriction.

In order to describe our main result about the homotopy type of the
classifying space $B \mathcal{C}_d$, we need some notation.  Let
$G(d,n)$ denote the Grassmannian of $d$-dimensional linear subspaces
of $\setR^{n+d}$.  There are two standard vector bundles, $U_{d,n}$
and $U_{d,n}^\perp$, over $G(d,n)$.  We are interested in the
$n$-dimensional one with total space
\begin{align*}
  U_{d,n}^\perp = \{(V,v) \in G(d,n) \times \setR^{d+n} \mid v \perp
  V\}.
\end{align*}
The Thom spaces (one-point compactifications) $\Th (U_{d,n}^\perp)$
define a spectrum $\mathit{MT}O(d)$ as $n$ varies\footnote{This
  convenient and flexible notation was suggested by Mike Hopkins.  In
  classical cobordism theory the standard notation for the Thom space
  of $U_{d,\infty} \to G(d,\infty)$ is $MO(d)$.  In that context,
  $O(d)$ is the structure group for normal bundles of manifolds.
  $O(d)$ is here the structure group for the tangent bundles of
  manifolds; hence the notation $MTO(d)$.}.
  The $(n+d)$th
space in the spectrum $\mathit{MT}O(d)$ is $\Th(U_{d,n}^\perp)$.  We
are primarily interested in the direct limit
\begin{align*}
  \Omega^{\infty-1}\mathit{MT}O(d) = \colim_{n \to \infty} \Omega^{n+d-1}
  \Th(U_{d,n}^\perp).
\end{align*}
$\mathit{MT}O(d)$ and $\Omega^{\infty-1}\mathit{MT}O(d)$ are described in more
detail in section~\ref{sec:spectrum-g_-d}.

Given a morphism $W \subseteq [a_0,a_1] \times \setR^{n+d-1}$, the
Pontrjagin-Thom collapse map onto a tubular neighborhood gives a map
from $[a_0, a_1]_+ \wedge S^{n+d-1}$ to the Thom space $\Th (\nu)$ of
the normal bundle of the embedding of $W$. Composing this with the
classifying map of $\nu$ yields a map
\begin{align*}
  [a_0, a_1]_+ \wedge S^{n+d-1} \to \Th(U_{d,n}^\perp),
\end{align*}
whose adjoint determines a path in $\Omega^{\infty-1} \mathit{MT}O(d)$ as $n \to
\infty$.  With more care, one gets a functor from $\mathcal{C}_d$ to
the category $\Path(\Omega^{\infty-1} \mathit{MT}O(d))$, whose objects are
points in $\Omega^{\infty-1} \mathit{MT}O(d)$ and whose morphisms are continuous
paths.

The classifying space of a path category is always homotopy equivalent
to the underlying space.  We therefore get a map
\begin{align}\label{eq:15}
  \alpha: B \mathcal{C}_d \to \Omega^{\infty-1} \mathit{MT}O(d)
\end{align}
(cf.\ \cite{MR1856399} for $d = 2$).

\begin{mthm}
  The map $\alpha: B \mathcal{C}_d \to \Omega^{\infty-1}
  \mathit{MT}O(d)$ is a weak homotopy equivalence.
\end{mthm}

For any category $C$, the set of components $\pi_0 BC$ can be
described as the quotient of the set $\pi_0 \ob(C)$ by the equivalence
relation generated by the morphisms.  For the category
$\mathcal{C}_d$, this gives that $\pi_0B\mathcal{C}_d$ is the group
$\Omega_{d-1}^O$ of cobordism classes of closed unoriented manifolds.
As explained in section~\ref{sec:spectrum-g_-d} below, the group of
components $\pi_0 \Omega^{\infty-1} \mathit{MT}O(d)$ is isomorphic to the
homotopy group $\pi_{d-1}MO$ of the Thom spectrum $MO$.  Thus the main
theorem can be seen as generalization of Thom's theorem: $\Omega_{d-1}^O
\cong \pi_{d-1}MO$.

More generally we also consider the cobordism category
$\mathcal{C}_\theta$ of manifolds with tangential structure, given by
a lifting of the classifying map for the tangent bundle over a
fibration $\theta: B \to G(d,\infty)$.  In this case, the right hand
side of~\eqref{eq:15} gets replaced by a spectrum
$\mathit{MT}(\theta)$ whose $(n+d)$th space is
$\Th(\theta^*U_{d,n}^\perp)$.  Chapter \ref{sec:tang-struct} defines
$\mathcal{C}_\theta$ and $\mathit{MT}(\theta)$ in more detail, and
proves the following version of the main theorem.
\begin{mthm2}
  There is a weak homotopy equivalence $\alpha^\theta: B
  \mathcal{C}_\theta \to \Omega^{\infty-1} \mathit{MT}(\theta)$.
\end{mthm2}

The simplest example of a tangential structure is that of an ordinary
orientation, leading to the category $\mathcal{C}_d^+$ of oriented
embedded cobordisms.  In this case, the target of $\alpha$ becomes the
oriented version $\Omega^{\infty-1} \mathit{MT}SO(d)$, which differs
from $\Omega^{\infty-1} \mathit{MT}O(d)$ only in that we start with
the Grassmannian $G^+(d,n)$ of oriented $d$-planes in $\setR^{n+d}$.
Another interesting special case leads to the category
$\mathcal{C}_d^{+}(X)$ of oriented manifolds with a continuous map to
a background space $X$.  In this case our result is a weak equivalence
\begin{align*}
  B\mathcal{C}_d^{+}(X) \simeq \Omega^{\infty-1} (\mathit{MT}SO(d) \wedge
  X_+).
\end{align*}
In particular, the homotopy groups $\pi_*B\mathcal{C}_d^{+}(X)$
becomes a generalized homology theory as a functor of the background
space $X$, with coefficients $\pi_* \Omega^{\infty-1}
\mathit{MT}SO(d)$.  The same works in the non-oriented situation.

We shall write $\G{d} = \mathit{MT}O(d)$ and $\G{d}^+ =
\mathit{MT}SO(d)$ for brevity, since we are mostly concerned with
these two cases.

For any topological category $\mathcal{C}$ and objects $x,y \in \ob
\mathcal{C}$, there is a continuous map
\begin{align*}
  \mathcal{C}(x,y) \to \Omega_{x,y} B \mathcal{C},
\end{align*}
from the space of morphisms in $\mathcal{C}$ from $x$ to $y$ to the
space $\Omega_{x,y} B \mathcal{C}$ of paths in $B\mathcal{C}$ from $x$
to $y$.  In the case of the oriented cobordism category we get for
every oriented $d$-manifold $W$ a map
\begin{align*}
  \sigma: B \Diff^+(W; \partial W) \to \Omega B\mathcal{C}_d^+
\end{align*}
into the loop space of $B\mathcal{C}_d^+$.  For $d=2$
and $W = W_{g,n}$ an oriented surface of genus $g$,
\begin{align*}
  B\Diff^{+}(W, \partial W) \simeq B \Gamma_{g,n},
\end{align*}
where $\Gamma_{g,n} = \pi_0 \Diff^+(W,\partial W)$ is the mapping
class group of $W$.  In this case, the composition
\begin{align*}
  B\Gamma_{\infty,n} \to \Omega_0 B \mathcal{C}_2^+ \xrightarrow{\simeq}
  \Omega^\infty_0 \mathit{MT}(2)^+
\end{align*}
induces an isomorphism in integral homology.  This is the generalized
Mumford conjecture, proved in \cite{arXiv1}.  We give a new proof of
this below, based on the above Main Theorem.


\section{The cobordism category and its sheaves}
\label{cha:cobord-categ-its}

\subsection{The cobordism category}
\label{sec:cobordism-category}

We fix the integer $d\geq 0$. The objects of the $d$-dimensional
cobordism category $\mathcal{C}_d$ are closed $(d-1)$-dimensional
smooth submanifolds of high-dimensional euclidean space; the morphisms
are $d$-dimensional embedded cobordisms with a collared boundary.

More precisely, an object of $\mathcal{C}_d$ is a pair $(M,a)$ with
$a\in \setR$, and such that $M$ is a closed $(d-1)$-dimensional
submanifold
\begin{align*}
  M \subseteq \setR^{d-1+\infty} \quad, \quad
  \setR^{d-1+\infty} = \colim_{n\to\infty} \setR^{d-1+n}
\end{align*}
A non-identity morphism from $(M_0,a_0)$ to $(M_1,a_1)$ is a triple
$(W,a_0,a_1)$ consisting of the numbers $a_0, a_1$, which must satisfy
$a_0 < a_1$, and a $d$-dimensional compact submanifold
\begin{align*}
  W \subseteq [a_0, a_1] \times \setR^{d-1+\infty},
\end{align*}
such that for some $\epsilon > 0$ we have
\begin{enumerate}[(i)]
\item $W \cap ([a_0,a_0 + \epsilon) \times \setR^{d-1+\infty}) = [a_0,
  a_0 + \epsilon) \times M_0$,
\item $W \cap ((a_1 -\epsilon ,a_1] \times \setR^{d-1+\infty}) = (a_1-
  \epsilon, a_1] \times M_1$,
\item $\partial W = W \cap (\{a_0, a_1\} \times \setR^{d-1 + \infty})$.
\end{enumerate}
Composition is union of subsets (of $\setR \times
\setR^{d-1+\infty}$):
\begin{align*}
  (W_1, a_0, a_1) \circ (W_2, a_1, a_2) = (W_1 \cup W_2, a_0, a_2).
\end{align*}
This defines $\mathcal{C}_d$ as a category of sets.  We describe its
topology.

Given a closed smooth $(d-1)$-manifold $M$, let $\Emb(M,
\setR^{d-1+n})$ denote the space of smooth embeddings, and write
\begin{align*}
  \Emb(M, \setR^{d-1+\infty}) &= \colim_{n\to\infty}
  \Emb(M, \setR^{d-1+n}).
\end{align*}
Composing an embedding with a diffeomorphism of $M$ gives a free
action of $\Diff(M)$ on the embedding space, and the orbit map
\begin{align*}
  \Emb(M, \setR^{d-1+\infty}) \to \Emb(M, \setR^{d-1+\infty})/\Diff(M)
\end{align*}
is a principal $\Diff(M)$ bundle in the sense of \cite{MR1688579}, if
$\Emb(M, \setR^{d-1+\infty})$ and $\Diff(M)$ are given Whitney
$C^\infty$ topology.

Let $E_\infty(M) = \Emb(M,\setR^{d-1+\infty})\times_{\Diff(M)} M$ and
let $B_\infty(M)$ be the orbit space
$\Emb(M,\setR^{d-1+\infty})/{\Diff(M)}$.  The associated fiber bundle
\begin{align}
  \label{eq:52}
  E_\infty(M) \to B_\infty(M)
\end{align}
has fiber $M$ and structure group $\Diff(M)$.  By Whitney's embedding
theorem $\Emb(M,\setR^{d-1+\infty})$ is contractible, so $B_\infty(M)
\simeq B \Diff(M)$.  In \cite{MR1471480} a convenient category of
infinite dimensional manifolds is described in which $\Diff(M)$ is a
Lie group and~(\ref{eq:52}) is a smooth fiber bundle.  The fiber
bundle~(\ref{eq:52}) comes with a natural embedding $E_\infty(M)
\subset B_\infty(M) \times \setR^{d-1+\infty}$.  With this structure,
it is universal.  More precisely, if $f: X \to B_\infty(M)$ is a
smooth map from a smooth manifold $X^k$, then the pullback
\begin{align*}
  f^*(E_\infty(M)) = \{ (x,v) \in X \times \setR^{d-1+\infty} \mid
  (f(x),v) \in E_\infty(M)\}
\end{align*}
is a smooth $(k+d)$-dimensional submanifold $E \subseteq X \times
\setR^{d-1+\infty}$ such that the projection $E \to X$ is a smooth
fiber bundle with fiber $M$.  Any such $E\subseteq X \times
\setR^{d-1+\infty}$ is induced by a unique smooth map $f: X \to
B_\infty(M)$.

Now the set of objects of $\mathcal{C}_d$ is
\begin{align}\label{eq:16}
  \ob\mathcal{C}_d \cong \setR \times \coprod_{M} B_\infty(M),
\end{align}
where $M$ varies over closed $(d-1)$-manifolds, one in each
diffeomorphism class.  We use this identification to topologize
$\ob\mathcal{C}_d$.

The set of morphisms in $\mathcal{C}_d$ is topologized in a similar
fashion.  Let $(W,h_0, h_1)$ be an abstract cobordism from $M_0$ to
$M_1$, i.e.\ a triple consisting of a smooth compact $d$-manifold $W$
and embeddings (``collars'')
\begin{align}\label{eq:6}
\begin{aligned}
  h_0: {} &[0,1) \times M_0 \to W\\
  h_1: {} &(0,1] \times M_1 \to W
\end{aligned}
\end{align}
such that $\partial W$ is the disjoint union of the two spaces
$h_\nu(\{\nu\} \times M_\nu)$, $\nu = 0,1$.  For $0 < \epsilon<
\frac12$, let $\Emb_\epsilon(W,[0,1] \times \setR^{d-1+n})$ be the
space of embeddings
\begin{align*}
  j: W \to [0,1]\times \setR^{d-1+n}
\end{align*}
for which there exists embeddings $j_\nu: M_\nu \to \setR^{d-1+n}$,
$\nu = 0,1$, such that
\begin{align*}
  j\circ h_0 (t_0,x_0) = (t_0,j_0(x_0)) \quad\text{and}\quad j\circ
  h_1 (t_1,x_1) = (t_1,j_1(x_1))
\end{align*}
for all $t_0 \in [0,\epsilon)$, $t_1 \in (1-\epsilon, 1]$, and $x_\nu
\in M_\nu$.  Let
\begin{align*}
  \Emb(W, [0,1] \times \setR^{d-1+\infty}) = \colim_{\substack{n \to
    \infty\\\epsilon \to 0}} \Emb_\epsilon(W, [0,1] \times \setR^{d-1+n}).
\end{align*}
Let $\Diff_\epsilon(W)$ denote the group of diffeomorphisms of $W$
that restrict to product diffeomorphisms on the $\epsilon$-collars,
and let $\Diff(W) = \colim_\epsilon \Diff_\epsilon(W)$.

As before, we get a principal $\Diff(W)$-bundle
\begin{align*}
  \Emb(W, [0,1]\times \setR^{d-1+\infty}) \to \Emb(W, [0,1]\times
  \setR^{d-1+\infty})/\Diff(W),
\end{align*}
and an associated fiber bundle
\begin{align*}
  E_\infty(W) \to B_\infty(W) = \Emb(W, [0,1]\times
  \setR^{d-1+\infty})/\Diff(W)
\end{align*}
with fiber $W$ and structure group $\Diff(W)$, satisfying a universal
property similar to the one for $E_\infty(M) \to B_\infty(M)$
described above.

Topologize $\mor\mathcal{C}_d$ by
\begin{align}
  \label{eq:7}
  \mor\mathcal{C}_d \cong \ob\mathcal{C}_d \amalg \coprod_{W}
  \setR^2_+ \times B_\infty(W),
\end{align}
where $\setR^2_+$ is the open half plane $a_0 < a_1$, and $W$ varies
over cobordisms $W = (W, h_0, h_1)$, one in each diffeomorphism class.

For $(a_0,a_1) \in \setR^2_+$, let $l: [0,1] \to [a_0, a_1]$ be the
affine map with $l(\nu) = a_\nu$, $\nu = 0,1$.  For an element $j \in
\Emb_\epsilon(W, [0,1] \times \setR^{d-1+\infty})$ we identify the
element $((a_0,a_1),[j]) \in \setR^2_+ \times B_\infty(W)$ with the
element $(a_0,a_1,E) \in \mor\mathcal{C}_d$, where $E$ is the image
\begin{align*}
  E = (l \circ j)(W) \subseteq [a_0, a_1] \times \setR^{d-1 + \infty}.
\end{align*}
Let us point out a slight abuse of notation: Strictly speaking, we
should include the collars $h_0$ and $h_1$ in the notation for the
$\Emb$ and $\Diff$ spaces.  Up to homotopy,
\begin{align}\label{diffsimdiff}
  \Diff(W) \xrightarrow{\simeq} \Diff(W,\{\partial_\mathrm{in} W\},
  \{\partial_\mathrm{out} W\})
\end{align}
is the group of diffeomorphisms of $W$ that restrict to
diffeomorphisms of the incoming and of the outgoing boundary of the
cobordism $W$.

Again, Whitney's embedding theorem implies that $B_\infty(W) \simeq
B\Diff(W)$.  With respect to this homotopy equivalence, composition in
$\mathcal{C}_d$ is induced by the morphism of topological groups
\begin{align*}
  \Diff(W_1) \times_{\Diff(M_1)} \Diff(W_2) \to \Diff(W),
\end{align*}
where $\partial_\mathrm{out} W_1 = M_1 = \partial_\mathrm{in} W_2$,
and $W = W_1 \cup_{M_1} W_2$.

\begin{remark}
  \label{rem:1}
  \begin{iblist}
  \item There is a reduced version $\widetilde{\mathcal{C}}_d$ where
    objects are embedded in $\{0\}\times \setR^{d-1+\infty}$ and
    morphisms in $[0, a_1]\times \setR^{d-1+\infty}$. The functor
    $\mathcal{C}_d\to\widetilde{\mathcal{C}}_d$ that maps a cobordism
    $W^d\subseteq [a_0, a_1]\times \setR^{d-1+\infty}$ into $W^d -
    a_0\in [0, a_1 - a_0] \times\setR^{d-1+\infty}$ induces a homotopy
    equivalence on classifying spaces. Indeed, the nerves are related
    by a pullback diagram
    \begin{align}\label{eq:21}
      \begin{aligned}
      \xymatrix{
        {N_k\mathcal{C}_d} \ar[r] \ar[d]&
        {N_k\widetilde{\mathcal{C}}_d} \ar[d]\\
        {N_k(\setR, \leq)} \ar[r] & {N_k(\setR_+, +)} 
      }
      \end{aligned}
    \end{align}
    where $(\setR, \leq)$ denotes $\setR$ as an ordered set and
    $(\setR_+,+)$ denotes $\setR_+ = \{0\} \amalg (0,\infty)$ as a
    monoid under addition.  The two vertical maps are fibrations, and
    the bottom horizontal map is a weak equivalence.  Therefore the
    functor $\mathcal{C}_d\to\widetilde{\mathcal{C}}_d$ induces a
    levelwise homotopy equivalence on nerves.
  \item In the previous remark it is crucial that $\setR$ be given its
    usual topology.  More precisely, let $\setR^\delta$ denote $\setR$
    with the discrete topology, and define $\mathcal{C}_d^\delta$ and
    $\widetilde{\mathcal{C}}_d^\delta$ using $\setR^\delta$ instead of
    $\setR$ in the homeomorphisms~\eqref{eq:16} and \eqref{eq:7}.
    Then the right hand vertical map in~\eqref{eq:21} defines a map
    $B\widetilde{\mathcal{C}}_d^\delta \to B(\setR_+^\delta,+)$ which is a
    split surjection.  By the group-completion theorem
    \cite{MR0402733}, $\pi_1 B(\setR_+^\delta,+) \cong \setR$, and this is a direct
    summand of $\pi_1B\widetilde{\mathcal{C}}_d^\delta$, so the main
    theorem fails for $\widetilde{\mathcal{C}}_d^\delta$.  We shall
    see later that $B\mathcal{C}_d^\delta \to B\mathcal{C}_d$ is a
    homotopy equivalence (cf.\
    Remark~\ref{remark:a-locally-constant}).
  \item There is a version $\mathcal{C}_d^+$ of $\mathcal{C}_d$ where
    one adds an orientation to the objects and morphisms in the usual
    way. For $d = 2$, the reduced version
    $\widetilde{\mathcal{C}}_d^+$ is the surface category
    $\mathcal{Y}$ of \cite[\S2]{MR1856399}.
  \end{iblist}
\end{remark}

\subsection{Recollection from \cite{arXiv1} on sheaves}
\label{sec:recollection-from-mw}

Let $\mathcal{X}$ denote the category of smooth (finite dimensional)
manifolds without boundary and smooth maps. We shall consider sheaves
on $\mathcal{X}$, that is, contravariant functors $\mathcal{F}$ on
$\mathcal{X}$ that satisfy the sheaf condition: for any open covering
$\mathcal{U} = \{ U_j \mid j\in J\}$ of an object $X$ in $\mathcal{X}$
and elements $s_j \in\mathcal{F}(U_j)$ with $s_j|U_i\cap U_j = s_i|
U_i\cap U_j$ there is a unique $s\in\mathcal{F}(X)$ that restricts to
$s_j$ for all $j$.  We have the Yoneda embedding of $\mathcal{X}$ into
the category $\ShX$ of sheaves on $\mathcal{X}$ that to $X \in
\mathcal{X}$ associates the representable sheaf $\tilde X =
C^\infty(-,X) \in \ShX$.

For the functors $\mathcal{F}$ we shall consider, $\mathcal{F}(X)$
consists of spaces over $X$ with extra properties.  In general the set
of spaces $E$ over $X$ is not a functor under pull-back ($(g\circ
f)^*(E) \not= f^*(g^*(E))$). But if $E\to X$ comes from subsets
$E\subseteq X\times U$ where $U$ is some ``universe'' then pull-backs with
respect to $f\colon X' \to X$ in $\mathcal{X}$, defined as
\begin{align*}
  f^*(E) &= \{ (x', u) \mid (f(x'), u) \in E\} \subseteq X'\times U, 
\end{align*}
is a functorial construction.

A set-valued sheaf $\mathcal{F}$ on $\mathcal{X}$ gives rise to a
\emph{representing space} $\abs{\mathcal{F}}$, constructed as the
topological realization of the following simplicial set. The
hyperplane (open or extended simplex)
\begin{align*}
  \Delta_e^\ell &= \{ (t_0, \dots, t_\ell) \in \setR^{\ell + 1} \mid
  \textstyle{\sum t_i} = 1\}
\end{align*}
is an object of $\mathcal{X}$, and
\begin{align*}
  [\ell] \longmapsto \mathcal{F}(\Delta_e^\ell)
\end{align*}
is a simplicial set. The space $\abs{\mathcal{F}}$ is its standard
topological realization. This is a representing space in the following
sense.
\begin{definition}
  \label{def:7}
  Two elements $s_0, s_1 \in \mathcal{F}(X)$ are concordant if there
  exists an $s\in \mathcal{F}(X\times\setR)$ which agrees with
  $\pr^*(s_0)$ in an open neighborhood of $X\times (-\infty, 0]$ and
  with $\pr^*(s_1)$ in an open neighborhood of $X\times [1, +\infty)$
  where $\pr\colon X\times\setR\to X$ is the projection.
\end{definition}

The set of concordance classes will be denoted by $\mathcal{F}[X]$.
The space $\abs{\mathcal{F}}$ above is a representing space in the
sense that $\mathcal{F}[X]$ is in bijective correspondence with the
set of homotopy classes of continuous maps from $X$ into
$\abs{\mathcal{F}}$:
\begin{align}
  \label{eq:32}
  \mathcal{F}[X] \cong [X, \abs{\mathcal{F}}]
\end{align}
by Proposition~A.1.1 of \cite{arXiv1}.  We describe the map.  For
$\tilde X = C^\infty (-,X)$, $[l]\mapsto \tilde X(\Delta^l_e)$ is the
(extended, smooth) total singular set of $X$, and satisfies that the
canonical map $|\tilde X| \to X$ is a homotopy equivalence
(\cite{MR0084138}).  An element $s \in \mathcal{F}(X)$ has an adjoint
$\tilde s: \tilde X \to \mathcal{F}$, inducing $|\tilde s|: |\tilde X|
\to |\mathcal{F}|$, and thus a well defined homotopy class of maps $X
\to |\mathcal{F}|$ which is easily seen to depend only on the
concordance class of $s$.

\begin{definition}
  \label{def:8}
  A map $\tau\colon \mathcal{F}_1\to\mathcal{F}_2$ is called a weak
  equivalence if the induced map from $\abs{\mathcal{F}_1}$ to
  $\abs{\mathcal{F}_2}$ induces an isomorphism on all homotopy groups.
\end{definition}

There is a convenient criteria for deciding if a map of sheaves is a
weak equivalence. This requires a relative version of
Definition~\ref{def:7}. Let $A\subseteq X$ be a closed subset of $X$,
and let $s\in \colim_U \mathcal{F}(U)$ where $U$ runs over open
neighborhoods of $A$. Let $\mathcal{F}(X, A; s) \subseteq
\mathcal{F}(X)$ be the subset of elements that agree with $s$ near
$A$.

\begin{definition}
  \label{def:9}
  Two elements $t_0, t_1\in\mathcal{F}(X, A; s)$ are concordant
  rel.\,$A$ if they are concordant by a concordance whose germ near
  $A$ is the constant concordance of $s$.  Let $\mathcal{F}_1[X, A;
  s]$ denote the set of concordance classes.
\end{definition}

\begin{criteria}
  \label{cri:1}
  A map $\tau\colon \mathcal{F}_1\to\mathcal{F}_2$ is a weak
  equivalence provided it induces a surjective map
  \begin{align*}
    \mathcal{F}_1[X, A; s] \to \mathcal{F}_2[X, A; \tau(s)]
  \end{align*}
  for all $(X, A, s)$ as above.
\end{criteria}

Let $x_0\in X$ and $s_0\in\mathcal{F}(\{x_0\})$. This gives a germ
$s_0\in\colim_U \mathcal{F}(U)$ with $U$ ranging over the open
neighborhoods of $x_0$. There is the following relative version of
\eqref{eq:32}, also proved in Appendix~A of \cite{arXiv1}: for every
$(X, A, s_0)$,
\begin{align*}
  \mathcal{F}[X, A; s_0] \cong [(X, A), (\abs{\mathcal{F}}, s_0)]
\end{align*}
In particular the homotopy groups $\pi_n(\abs{\mathcal{F}}, s_0)$ are
equal to the relative concordance classes $\mathcal{F}[S^n, x_0;
s_0]$. By Whitehead's theorem $\tau\colon
\mathcal{F}_1\to\mathcal{F}_2$ is a weak equivalence if and only if
\begin{align*}
  \mathcal{F}_1[S^n, x_0; s_0] \xrightarrow{\ \cong\ } \mathcal{F}_2[S^n,
  x_0; \tau(s_0)]
\end{align*}
is an equivalence for all basepoints $x_0$ and all
$s_0\in\mathcal{F}_1(x_0)$. This is sometimes a more convenient
formulation than Criteria~\ref{cri:1} above.

Actually, for the concrete sheaves we consider in this paper the
representing spaces are ``simple'' in the sense of homotopy theory,
and in this situation the base point $s_0\in\mathcal{F}(*)$ is
irrelevant: a map $\tau\colon\mathcal{F}_1\to\mathcal{F}_2$ is a weak
equivalence if and only if it induces a bijection $\mathcal{F}_1[X]
\to\mathcal{F}_2[X]$ for all $X \in \mathcal{X}$. In fact, it suffices
to check this when $X$ is a sphere.

\subsection{A sheaf model for the cobordism category}
\label{sec:sheaf-model}

We apply the above to give a sheaf model of the cobordism category
$\mathcal{C}_d$. First some notation. For functions $a_0, a_1\colon
X\to \setR$ with $a_0(x) \leq a_1(x)$ at all $x\in X$, we write
\begin{align*}
  X\times (a_0, a_1) &= \{ (x, u)\in X\times\setR \mid a_0(x) < u <
  a_1(x)\}\\
  X\times [a_0, a_1] &= \{ (x, u) \in X\times \setR \mid a_0(x) \leq u
  \leq a_1(x)\}.
\end{align*}
Given a submersion $\pi\colon W \to X$ of smooth manifolds (without
boundary) and smooth maps
\begin{align*}
  f\colon W\to \setR, &\quad a\colon X\to\setR,
\end{align*}
we say that $f$ is \emph{fiberwise transverse} to $a$ if the
restriction $f_x$ of $f$ to $W_x = \pi\inv(x)$ is transversal to
$a(x)$ for every $x\in X$, or equivalently if the graph $X\times\{a\}$
consists of regular values for $(\pi,f):E\to X\times \setR$.
 In this case
\begin{align*}
  M = (f - a\pi)^{-1}(0) = \{ z\in W \mid f(z) = a(\pi(z))\}
\end{align*}
is a codimension one submanifold of $W$, and the restriction
$\pi\colon M\to X$ is still a submersion.

For $X \in \mathcal{X}$ and smooth real functions
\begin{align*}
  a_0 \leq a_1 \colon X\to \setR, &\quad \epsilon\colon X\to (0,
  \infty)
\end{align*}
we shall consider submanifolds
\begin{align*}
  W \subseteq X \times (a_0 - \epsilon, a_1 + \epsilon) \times
  \setR^{d-1+\infty}.
\end{align*}
The three projections will be denoted
\begin{align*}
  \pi\colon W\to X, \quad f\colon W\to \setR, \quad j\colon W\to
  \setR^{d-1+\infty}
\end{align*}
unless otherwise specified.

\begin{definition}
  \label{def:3}
  For $X \in \mathcal{X}$ and smooth real functions $a_0 \leq a_1$
  and $\epsilon$ as above, the set $C_d^\pitchfork(X; a_0, a_1, \epsilon)$
  consists of all submanifolds
  \begin{align*}
    W \subseteq X \times (a_0-\epsilon,
    a_1+\epsilon)\times\setR^{d-1+\infty}
  \end{align*}
  which satisfies the following conditions:
  \begin{enumerate}[(i)]
  \item $\pi\colon W\to X$ is a submersion with $d$-dimensional
    fibers,
  \item $(\pi, f)\colon W\to X\times (a_0-\epsilon, a_1+\epsilon)$ is proper,
  \item The restriction of $(\pi,f)$ to $(\pi,f)^{-1}(X \times (a_\nu-
    \epsilon, a_\nu + \epsilon))$ is a submersion for $\nu = 0,1$.
  \end{enumerate}
\end{definition}

The three conditions imply that $\pi\colon W\to X$ is a smooth fiber
bundle rather than just a submersion. Indeed for each $\nu = 0,1$,
restricting $(\pi,f)$ gives a map
\begin{align*}
  (\pi,f)^{-1}(X \times (a_\nu- \epsilon, a_\nu + \epsilon)) \to X
  \times (a_\nu- \epsilon, a_\nu + \epsilon)
\end{align*}
which is a proper submersion, and hence a smooth fiber bundle by
Ehresmann's fibration lemma, cf.\ \cite[p.~84]{MR674117}. Similarly
the restriction of $\pi$ to
\begin{align*}
  W[a_0, a_1] &= W\cap X\times [a_0, a_1] \times \setR^{d-1+\infty}
\end{align*}
is a smooth fiber bundle with boundary. The result for $\pi\colon W\to
X$ follows by gluing the collars.

We remove the dependence on $\epsilon$ and define
\begin{align*}
  C_d^\pitchfork(X; a_0, a_1) &= \colim_{\epsilon\to 0}
  C_d^\pitchfork(X; a_0, a_1, \epsilon).
\end{align*}
\begin{definition}
  \label{def:2}
  For $X\in\mathcal{X}$, let
  \begin{align*}
    C_d^\pitchfork(X) &= \coprod C_d^\pitchfork(X; a_0, a_1).
  \end{align*}
  The disjoint union varies over the (uncountable) set consisting of
  pairs of smooth functions with $a_0 \leq a_1$ and such that $\{x\mid
  a_0(x) = a_1(x)\}$ is open (hence a union of connected components of
  $X$).  This defines a sheaf $C_d^\pitchfork$.
\end{definition}

Taking union of embedded manifolds gives a partially defined map
\begin{align*}
  C_d^\pitchfork(X; a_0, a_1) \times C_d^\pitchfork(X; a_1, a_2) \to
  C_d^\pitchfork(X; a_0, a_2)
\end{align*}
and defines a category structure on $C_d^\pitchfork(X)$ with the
objects (or identity morphisms) corresponding to $a_0 = a_1$.

A smooth map $\phi\colon Y\to X$ induces a map of categories $\phi^*
\colon C_d^\pitchfork(X) \to C_d^\pitchfork(Y)$ by the pull-back
construction of \S~\ref{sec:recollection-from-mw}: For
\begin{align*}
  W \subseteq X\times (a_0-\epsilon, a_1+\epsilon)\times\setR^{d-1+\infty},
\end{align*}
$\phi^* W = \{ (y, u, r) \mid (\phi(y), u, r) \in W\}$ is an element
of $C_d^\pitchfork(Y; a_0\phi, a_1\phi, \epsilon\phi)$. This gives a
$\CAT$-valued sheaf
\begin{align*}
  C_d^\pitchfork \colon\mathcal{X}\to \CAT,
\end{align*}
where $\CAT$ is the category of small categories.

An object of $C_d^\pitchfork(\point)$ is represented by a $d$-manifold
$W\subseteq (a-\epsilon,a+\epsilon) \times \setR^{d-1+\infty}$ such
that $f: W \to (a-\epsilon, a+\epsilon)$ is a proper submersion.  Thus
$M = f^{-1}(0)\subseteq \setR^{d-1+\infty}$ is a closed
$(d-1)$-manifold.  Only the \emph{germ} of $W$ near $M$ is
well-defined.  As an abstract manifold, $W$ is diffeomorphic to $M
\times (a-\epsilon,a+\epsilon)$, but the embedding into
$(a-\epsilon,a+\epsilon)\times \setR^{d-1+\infty}$ need not be the
product embedding.  Hence the germ of $W$ near $M$ carries slightly
more information than just the submanifold $M\subseteq \{a\} \times
\setR^{d-1+\infty}$.  This motivates
\begin{definition}\label{def:1}
  Let $C_d(X; a_0, a_1, \epsilon)\subseteq C_d^\pitchfork(X; a_0, a_1,
  \epsilon)$ be the subset satisfying the further condition
  \begin{enumerate}[(i)]
    \setcounter{enumi}{3}
  \item For $x \in X$ and $\nu = 0,1$, let $J_\nu$ be the interval
    $((a_\nu - \epsilon)(x), (a_\nu + \epsilon)(x))\subseteq \setR$,
    and let $V_\nu = (\pi,f)^{-1}(\{x\}\times J_\nu)\subseteq \{x\}
    \times J_\nu \times \setR^{d-1+\infty}$.  Then
    \begin{align*}
      V_\nu = \{x\} \times J_\nu \times M
    \end{align*}
    for some $(d-1)$-dimensional submanifold $M \subseteq
    \setR^{d-1+\infty}$.
  \end{enumerate}
  Define $C_d(X; a_0, a_1)\subseteq C^\pitchfork_d(X; a_0, a_1)$ and
  $C_d(X)\subseteq C_d^\pitchfork(X)$ similarly.
\end{definition}

It is easy to see that $C_d(X)$ is a full subcategory of
$C_d^\pitchfork(X)$ and that
\begin{align*}
  C_d\colon \mathcal{X} \to \CAT
\end{align*}
is a sheaf of categories,  isomorphic to the sheaf
$C^\infty(-,\mathcal{C}_d)$, where $\mathcal{C}_d$ is equipped with
the (infinite dimensional) smooth structure described in
section~\ref{sec:cobordism-category}.  In particular we get a
continuous functor
\begin{align*}
  \eta: |C_d| \to \mathcal{C}_d.
\end{align*}

\begin{proposition}\label{prop:B-eta}
  $B\eta: B|C_d| \to B\mathcal{C}_d$ is a weak homotopy equivalence.
\end{proposition}
\begin{proof}
  The space $N_k|C_d|$ is the realization of the simplicial set
  \begin{align*}
    [l] \mapsto N_k C_d(\Delta^l_e) = C^\infty(\Delta^l_e,
    N_k\mathcal{C_d}).
  \end{align*}
  A theorem from \cite{MR0084138} asserts that the realization of the
  singular simplicial set of any space $Y$ is weakly homotopy
  equivalent to $Y$ itself.  This is also the case if one uses the
  extended simplices $\Delta_e^k$ to define the singular simplicial
  set, and for manifolds it is also true if we use smooth maps.  This
  proves that the map
  \begin{align*}
    N_k \eta: N_k |C_d| \to N_k \mathcal{C}_d
  \end{align*}
  is a weak homotopy equivalence for all $k$, and hence that $B\eta$
  is a weak homotopy equivalence.
\end{proof}

\subsection{Cocycle sheaves}
\label{sec:cocycle-sheaves}

We review the construction from \cite[\S4.1]{arXiv1} of a model for
the classifying space construction at the sheaf level.

Let $\mathcal{F}$ be any $\CAT$-valued sheaf on $\mathcal{X}$. There
is an associated set valued sheaf $\beta\mathcal{F}$.  Choose, once
and for all, an uncountable set $J$.  An element of
$\beta\mathcal{F}(X)$ is a pair $(\mathcal{U}, \Phi)$ where
$\mathcal{U} = \{ U_j \mid j\in J\}$ is a locally finite open cover of
$X$, indexed by $J$, and $\Phi$ a certain collection of morphisms. In
detail: given a non-empty finite subset $R\subseteq J$, let $U_R$ be the
intersection of the $U_j$'s for $j\in R$. Then $\Phi$ is a collection
$\phi_{RS}\in N_1\mathcal{F}(U_S)$ indexed by pairs $R\subseteq S$ of
non-empty finite subsets of $J$, subject to the conditions
\begin{enumerate}[(i)]
\item $\phi_{RR} = \id_{c_R}$ for an object $c_R \in
  N_0\mathcal{F}(U_R)$,
\item For each non-empty finite $R\subseteq S$, $\phi_{RS}$ is a
  morphism from $c_S$ to $c_R | U_S$,
\item For all triples $R \subseteq S \subseteq T$ of finite non-empty
  subsets of $J$, we have
  \begin{align}
    \label{eq:38}
    \phi_{RT} = (\phi_{RS}|U_T) \circ \phi_{ST}.
  \end{align}
\end{enumerate}
Theorem~4.1.2 of \cite{arXiv1} asserts a weak homotopy equivalence
\begin{align}
  \label{eq:39}
  \abs{\beta\mathcal{F}} \simeq B\abs{\mathcal{F}}.
\end{align}
\begin{remark}
  \label{rem:4}
  In the case $\mathcal{F}(X) = \Map(X, \mathcal{C})$ for some
  topological category $\mathcal{C}$ the construction
  $\beta\mathcal{F}$ takes the following form. Let $X_{\mathcal{U}}$
  be the topological category from \cite{MR0232393}:
  \begin{align*}
    \ob X_{\mathcal{U}} &= \coprod_{R} U_R & \mor X_{\mathcal{U}} &=
    \coprod_{R\subseteq S} U_S,
  \end{align*}
  i.e.\ $X_{\mathcal{U}}$ is the topological poset of pairs $(R,x)$,
  where $R \subseteq J$ is a finite non-empty subset and $x \in U_R$.
  If $R \subseteq S$ and $x = y$, then there is precisely one morphism
  $(S,x) \to (R,y)$, otherwise there is none.

  Then \eqref{eq:38} amounts to a continuous functor $\Phi\colon
  X_{\mathcal{U}} \to \mathcal{C}$.  In general, (\ref{eq:38}) amounts
  to a functor $\tilde X_{\mathcal{U}} \to \mathcal{F}$, where $\tilde
  X_{\mathcal{U}} = C^\infty(-, X_{\mathcal{U}})$ is the
  (representable) sheaf of posets associated to $X_{\mathcal{U}}$.

  A partition of unity $\{ \lambda_j \mid j\in J\}$ subordinate to
  $\mathcal{U}$ defines a map from $X$ to $BX_{\mathcal{U}}$ and
  $\Phi$ a map from $BX_{\mathcal{U}}$ to $B\mathcal{C}$. This induces
  a map
  \begin{align*}
    \beta\mathcal{F}[X] \to [X, B\mathcal{C}]
  \end{align*}
  and \eqref{eq:39} asserts that this is a bijection for all $X$.
\end{remark}


\section{The Thom spectra and their sheaves}
\label{cha:thom-spectra-their}

\subsection{The spectrum $\G{d}$ and its infinite loop space}
\label{sec:spectrum-g_-d}

We write $G(d, n)$ for the Grassmann manifold of $d$-dimensional
linear subspaces of $\setR^{d+n}$ and $G^+(d, n)$ for the double cover
of $G(d, n)$ where the subspace is equipped with an orientation.

There are two distinguished vector bundles over $G(d, n)$, the
tautological $d$-dimensional vector bundle $U_{d, n}$ consisting of
pairs of a $d$-plane and a vector in that plane, and its orthogonal
complement, the $n$-dimensional vector bundle $\udnb$. The
direct sum $U_{d, n}\oplus \udnb$ is the product bundle $G(d,
n)\times \setR^{d+n}$.

The Thom spaces (one point compactifications) $\Th(\udnb)$ form the
spectrum $\G{d}$ as $n$ varies. Indeed, since $U_{d, n+1}^\bot$
restricts over $G(d, n)$ to the direct sum of $\udnb$ and a trivial
line, there is an induced map
\begin{align}
  \label{eq:17}
  S^1 \wedge \Th(\udnb) \to \Th(U_{d, n+1}^\bot).
\end{align}
The $(n+d)$th space of the spectrum $\G{d}$ is $\Th(U_{d, n}^\bot)$,
and~\eqref{eq:17} provides the structure maps. The associated infinite
loop space is therefore
\begin{align*}
  \Omega^\infty\G{d} = \colim_{n\to \infty} \Omega^{n+d} \Th(\udnb),
\end{align*}
where the maps in the colimit
\begin{align*}
  \Omega^{n+d}\Th(\udnb) \to \Omega^{n+d+1} \Th(U_{d, n+1}^\bot)
\end{align*}
are the $(n+d)$-fold loops of the adjoints of~\eqref{eq:17}.

There is a corresponding oriented version $\G{d}^+$ where one uses the
Thom spaces of pull-backs $\theta^*\udnb$, $\theta\colon G^+(d, n) \to
G(d, n)$. The spectrum $\G{d}^+$ maps to $\G{d}$ and induces
\begin{align*}
  \Omega^\infty\G{d}^+ \to \Omega^\infty\G{d}.
\end{align*}

\begin{proposition}
  \label{prop:2}
  There are homotopy fibration sequences
  \begin{align*}
    \Omega^\infty\G{d} &\xrightarrow{}
    \Omega^\infty\Sigma^\infty(BO(d)_+) \xrightarrow{\ \partial\ }
    \Omega^\infty\G{d-1},\\
    \Omega^\infty\G{d}^+ &\xrightarrow{}
    \Omega^\infty\Sigma^\infty(BSO(d)_+) \xrightarrow{\ \partial\ }
    \Omega^\infty\G{d-1}^+.
  \end{align*}
\end{proposition}
\begin{proof}
  For any two vector bundles $E$ and $F$ over the same base $B$ there
  is a cofiber sequence
  \begin{align}
    \label{eq:20}
    \Th(p^*E) \to \Th(E) \to \Th(E\oplus F)
  \end{align}
  where $p\colon S(F)\to X$ is the bundle projection of the sphere
  bundles.

  Apply this to $X = G(d, n)$, $E = \udnb$, $F = U_{d, n}$. The sphere
  bundle is
  \begin{align*}
    S(U_{d,n}) = O(n+d)/O(n)\times O(d-1).
  \end{align*}
  Since $G(d-1, n) = O(n+d-1)/O(n)\times O(d-1)$, the natural map
  $G(d-1, n)\to S(U_{d, n})$ is $(n+d-2)$-connected. The bundle $p^*
  \udnb$ over $S(U_{d, n})$ restricts to $U_{d-1, n}^\bot$ over
  $G(d-1, n)$, so
  \begin{align*}
    \Th(U_{d-1, n}^\bot) \to \Th(p^*\udnb)
  \end{align*}
  is $(2n+d-2)$-connected. The right-hand term in~\eqref{eq:20} is
  $G(d, n)_+ \wedge S^{n+d}$, and the map $G(d, n)\to BO(d)$ is
  $(n-1)$-connected ($BO(d) = G(d, \infty)$).

  The cofiber sequence~\eqref{eq:20} gives a cofiber sequence of
  spectra
  \begin{align}\label{eq:54}
    \Sigma\inv\G{d-1} \to \G{d} \to
    \Sigma^\infty(BO(d)_+) \to \G{d-1}
  \end{align}
  and an associated homotopy fibration sequence
  \begin{align*}
    \Omega^\infty\G{d} \to \Omega^\infty\Sigma^\infty(BO(d)_+) \to
    \Omega^\infty\G{d-1}
  \end{align*}
  of infinite loop spaces. The oriented case is completely similar.
\end{proof}

\begin{remark}
  \label{rem:3}
  For $d = 1$, the sequences in Proposition~\ref{prop:2} are
  \begin{align*}
    \Omega^\infty\G{1} &\xrightarrow{}
    \Omega^\infty\Sigma^\infty(\RP_+^\infty) \xrightarrow{\ \partial\ }
    \Omega^\infty\Sigma^\infty\\
    \Omega^\infty\G{1}^+ &\xrightarrow{} \Omega^\infty\Sigma^\infty
    \xrightarrow{\ \partial\ } \Omega^\infty\Sigma^\infty\times\Omega^\infty\Sigma^\infty.
  \end{align*}
  In the first sequence, $\partial$ is the stable transfer associated
  with the universal double covering space.  In the oriented case,
  $\partial$ is the diagonal. Thus
  \begin{align*}
    \Omega^\infty\G{1} = \Omega^\infty\RP_{-1}^\infty, \quad
    \Omega^\infty\G{1}^+ = \Omega(\Omega^\infty\Sigma^\infty).
  \end{align*}
  The oriented Grassmannian $G^+(2,\infty)$ is homotopy equivalent to
  $\setC P^\infty$, and the space $\Omega^\infty\G{2}^+$ is homotopy
  equivalent to the space $\Omega^\infty\CP_{-1}^\infty$, in the
  notation from \cite{arXiv1}.
\end{remark}

The cofiber sequence~\eqref{eq:54} defines a direct system of spectra
\begin{align}\label{eq:53}
  \G{0} \to \Sigma \G{1} \to \dots \to \Sigma^{d-1}\G{d-1} \to
  \Sigma^d \G{d} \to \cdots
\end{align}
whose direct limit
we could denote $MTO$.  In fact it is homotopy equivalent to the
universal Thom spectrum usually denoted $MO$ in the following way.
There is a homeomorphism $G(d,n) \to G(n,d)$ covered by a bundle
isomorphism $U^\perp_{d,n} \to U_{n,d}$.  Thus we have a maps
\begin{align}\label{eq:35}
  \Th(U^\perp_{d,n}) \xrightarrow{\cong} \Th(U_{n,d}) \to
  \Th(U_{n,\infty}).
\end{align}
The spaces $\Th(U^\perp_{d,n})$ and $\Th(U_{n,\infty})$ are the $n$th
spaces of the spectra $\Sigma^d\G{d}$ and $MO$, respectively, and the
map~\eqref{eq:35} induces a map of spectra $\Sigma^d \G{d} \to MO$.
$\Th(U_{n,\infty})$ can be built from $\Th(U_{n,d})$ by attaching
cells of dimension greater than $n+d$, so the resulting map $\Sigma^d
\G{d} \to MO$ induces an isomorphism in $\pi_k$ for $k < d$ and a
surjection for $k = d$.

The homotopy groups of $MO$ form the unoriented bordism ring
\begin{equation*}
  \pi_{d-1} MO = MO_{d-1}(\point) = \Omega^O_{d-1}.
\end{equation*}
The direct system~\eqref{eq:53} can be thought of as a filtration of
$MO$, with filtration quotients $\Sigma^d BO(d)_+$.  In particular,
the maps in the direct system induce an isomorphism
\begin{align*}
  \pi_{-1} \G{d} = \pi_{d-1} \Sigma^d \G{d}
  \stackrel{\isom}{\longrightarrow} \pi_{d-1} MO = \Omega^O_{d-1},
\end{align*}
and an exact sequence
\begin{align}\label{eq:55}
  \pi_0\G{d+1} \stackrel{\chi}{\longrightarrow} \setZ
  \stackrel{S^d}{\longrightarrow} \pi_0\G{d} \to \Omega^{O}_d \to 0.
\end{align}
The map $\chi: \pi_0\G{d+1} \to \setZ$ corresponds under the
homotopy equivalence of our main theorem to the map that to a closed
$(d+1)$-manifold $W$, thought of as an endomorphism in
$\mathcal{C}_{d+1}$ of the empty $d$-manifold, associates the Euler
characteristic $\chi(W) \in \setZ$.  The map $S^d: \setZ \to
\pi_0\G{d}$ corresponds to the $d$-sphere $S^d$, thought of as an
endomorphism in $\mathcal{C}_d$ of the empty $(d-1)$-manifold.  For
odd $d$, $\chi$ is surjective ($\chi(\RP^{d+1}) = 1$), so the
sequence~\eqref{eq:55} defines an isomorphism $\pi_0\G{d} \isom
\Omega^O_{d}$.  On the other hand $\chi=0$ for even $d$ by Poincar\'e
duality, so the sequence~\eqref{eq:55} works out to be
\begin{align*}
  0 \to \setZ \stackrel{S^d}{\longrightarrow} \pi_0\G{d} \to
  \Omega^O_d  \to 0.
\end{align*}

\subsection{Using Phillips' submersion theorem}
\label{sec:using-phill-subm}

We give a sheaf model for the space $\Omega^{\infty-1}\G{d}$.

\begin{definition}
  \label{def:6}
  For a natural number $n > 0$ and $X\in\mathcal{X}$, an element of
  $D_d(X; n)$ is a submanifold
  \begin{align*}
    W \subseteq X\times\setR\times\setR^{d-1+n},
  \end{align*}
  with projections $\pi$, $f$, and $j$, respectively, such that
  \begin{enumerate}[(i)]
  \item\label{item:1} $\pi\colon W\to X$ is a submersion with $d$-dimensional
    fibers.
  \item\label{item:2} $(\pi, f)\colon W \to X\times\setR$ is proper.
  \end{enumerate}
  This defines a set valued sheaf $D_d(-;n) \in \ShX$.  Let $D_d$ be
  the colimit (in $\ShX$) of $D_d(-;n)$ as $n\to\infty$.  Explicitly,
  $D_d(X)$ is the set of submanifolds $W \subseteq
  X\times\setR\times\setR^{d-1+\infty}$ satisfying (\ref{item:1}) and
  (\ref{item:2}) above, and such that for each compact $K \subseteq X$
  there exists an $n$ with $\pi^{-1}(K)\subseteq K \times \setR \times
  \setR^{d-1+n}$.
\end{definition}
We will prove the following theorem by constructing a natural
bijection $[X, \Omega^{\infty-1}\G{d}] \cong D_d[X]$.
\begin{theorem}
  \label{thm:4}
  There is a weak homotopy equivalence
  \begin{align*}
    \abs{D_d} \xrightarrow{\ \simeq\ } \Omega^{\infty-1}\G{d}.
  \end{align*}
\end{theorem}

Given $W\subseteq X\times\setR\times\setR^{d-1+n}$ with $n$-dimensional
normal bundle $N \to W$, there is a vector bundle map
\begin{align}
  \label{eq:22}
  \begin{aligned}
    \xymatrix{N \ar[r]^-{\hat\gamma} \ar[d] & U_{d, n}^\bot \ar[d]\\
        W\ar[r]_-{\gamma} & G(d, n).}
  \end{aligned}
\end{align}
Write $W_x$ for the intersection $W_x = W\cap \{x\}\times
\setR\times\setR^{d-1+n}$.  Then $\gamma(z) = T_z(W_{\pi(z)})$,
considered as a subspace of $\setR^{d+n}$. The normal fiber $N_z$ of
$W$ in $X\times\setR\times\setR^{d-1+n}$ is the normal fiber of $W_x$
in $\setR^{d+n}$, so is equal to $\gamma(z)^\bot$; this defines
$\hat\gamma$ in \eqref{eq:22}.

Next we pick a regular value for $f\colon W\rightarrow \setR$, say
$0\in\setR$, and let $M = f\inv(0)$.  Then the normal bundle $N$ of $W
\subseteq X \times \setR \times \setR^{d-1+n}$ restricts to the normal
bundle of $M\subset X\times\setR^{d-1+n}$.  Choose a tubular
neighborhood of $M$ in $X\times\setR^{d-1+n}$, and let
\begin{align*}
  e\colon {N|M} \to X\times \setR^{d-1+n}
\end{align*}
be the associated embedding (\cite[\S12]{MR674117}). The induced map of
one-point compactifications, composed with~\eqref{eq:22}, gives a map
\begin{align}
  \label{eq:23}
  g\colon X_+ \wedge S^{d-1+n} \to \Th(\udnb)
\end{align}
whose homotopy class is independent of the choices made (when $n\gg
d$). Its adjoint is a well-defined homotopy class of maps from $X$ to
$\Omega^{\infty-1}\G{d}$. This defines
\begin{align*}
  \rho\colon D_d[X] \to [X, \Omega^{\infty-1}\G{d}].
\end{align*}
We now construct an inverse to $\rho$ using transversality and
Phillips' submersion theorem.  We give the argument only in the case
where $X$ is compact.  Any map~\eqref{eq:23} is homotopic to a map
that is transversal to the zero section, and
\begin{align*}
  M= g\inv(G(d, n))  \subseteq X \times \setR^{d-1+n}
\end{align*}
is a submanifold.  The projection $\pi_0\colon M\to X$ is proper, and
the normal bundle is $N = g^*(\udnb)$. Define $T^\pi M = g^*(U_{d,
  n})$ so that
\begin{align*}
  N \oplus T^\pi M = M\times \setR^{n+d}.
\end{align*}
Combined with the bundle information of the embedding of $M$ in
$X\times\setR^{d-1+n}$ this yields an isomorphism of vector bundles
over $M$
\begin{align}
  \label{eq:25}
  TM \times \setR^{n+d} \xrightarrow{\ \cong\ } (\pi_0^* TX \oplus
  T^\pi M) \times \setR^{d-1+n}.
\end{align}

By standard obstruction theory (cf.\ \cite{arXiv1}, Lemma~3.2.3) there
is an isomorphism (unique up to concordance)
\begin{align*}
  \hat\pi_0 \colon TM\times \setR \xrightarrow{\ \cong\ } \pi_0^* TX\oplus
  T^\pi M
\end{align*}
that induces~\eqref{eq:25}. Set $W = M\times \setR$, $\pi_1 = \pi_0
\circ \pr_M$ and $T^\pi W = \pr_M^* T^\pi M$. Then
\begin{align}
  \label{eq:26}
  TW \xrightarrow{\ \cong\ } \pi_1^* TX \oplus T^\pi W, 
\end{align}
and since $W$ has no closed components we are in a position to apply
the submersion theorem. Indeed,~\eqref{eq:26} gives a bundle
epimorphism $\hat\pi_1\colon TW\to TX$ over $\pi_1\colon W\to X$. By
Phillips' theorem, there is a homotopy $(\pi_t,\hat \pi_t), t \in
[1,2]$ through bundle epimorphisms, from $(\pi_1, \hat\pi_1)$ to a
pair $(\pi_2, d\pi_2)$, i.e. to a submersion $\pi_2$.  Let $f: W \to
\setR$ be the projection.  Then $(\pi_2,f): W \to X \times \setR$ is
proper since we have assumed that $X$ is compact.  For $n \gg d$ we
get an embedding $W \subset X \times \setR \times \setR^{d-1+n}$ which
lifts $(\pi_2,f)$.

If $n\gg d$ the original embedding $W\subset
X\times\setR\times\setR^{d-1+\infty}$ is isotopic to an embedding
where the projection onto $X$ is the submersion $\pi$ and with $(\pi,
f)$ proper. (This is direct from \cite{MR0208611} when $X$ is compact;
and in general a slight extension.) We have constructed
\begin{align}
  \label{eq:27}
  \sigma\colon [X, \Omega^{\infty-1}\G{d}] \to D_d[X].
\end{align}
\begin{proposition}
  \label{prop:3}
  The maps $\sigma$ and $\rho$ are inverse bijections.
\end{proposition}
\begin{proof}
  By construction $\rho\circ\sigma = \id$. The other composite
  $\sigma\circ\rho = \id$ uses that an element $W\in D_d(X)$ is
  concordant to one where $W$ is replaced by $M\times\setR$ and $f$ by
  the projection; $M$ is the inverse image of a regular value of $f$.
  The concordance is given in Lemma~2.5.2 of \cite{arXiv1}.
\end{proof}

\begin{remark}
  One can define $\sigma$ also for non-compact $X$, but it requires a
  slight extension of \cite{MR0208611} to see that $(\pi_2,f): W \to X
  \times \setR$ can be taken to be proper.  The proof of
  Theorem~\ref{thm:4} above only uses~\eqref{eq:27} for compact $X$,
  in fact for $X$ a sphere.
\end{remark}


\section{Proof of the main theorem}
\label{sec:proof-of-main-theorem}

The proof uses an auxiliary sheaf of categories $D_d^\pitchfork$ and a
zig-zag of functors
\begin{align*}
  D_d \stackrel{\alpha}{\longleftarrow} \ddp
  \stackrel{\gamma}{\longrightarrow} \cdp
  \stackrel{\delta}{\longleftarrow} C_d
\end{align*}
The sheaf $C_d$ is the cobordism category sheaf, defined in
section~\ref{sec:sheaf-model} above, and $\cdp$ is the slightly larger
sheaf, defined in the same section.  The sheaf $D_d$ is, by
Theorem~\ref{thm:4}, a sheaf model of $\Omega^{\infty-1}\G{d}$.  We
regard $D_d$ as a sheaf of categories with only identity morphisms.
To prove the main theorem it will suffice to prove that $\alpha$,
$\gamma$ and $\delta$ all induce weak equivalences.

\begin{definition}
  \label{def:4}
  Let $\ddp(X)$ denote the set of pairs $(W,a)$ such that
  \begin{enumerate}[(i)]
  \item $W \in D_d(X)$,
  \item $a: X \to \setR$ is smooth,
  \item $f: W \to \setR$ is fiberwise transverse to $a$.
  \end{enumerate}
  
  Thus, $\ddp$ is a subsheaf of $D_d \times \tilde\setR$, where
  $\tilde \setR$ is the representable sheaf $C^\infty(-,\setR)$.  It
  is also a sheaf of posets, where $(W,a) \leq (W',a')$ when $W = W'$,
  $a\leq a'$ and $(a' - a)^{-1}(0)\subseteq X$ is open.
\end{definition}

Recall from section~\ref{sec:sheaf-model} that $f: W \to \setR$ is
fiberwise transverse to $a: X \to \setR$ if $f_x \colon W_x \to \setR$
is transverse to $a(x)\in \setR$ for all $x \in X$.
By properness of $(\pi,f)$, there will exist a smooth map
$\epsilon: X \to (0,\infty)$, such that the restriction of $(\pi,f)$
to the open subset
\begin{align*}
  W_\epsilon = (\pi,f)^{-1}(X \times (a-\epsilon, a+\epsilon)),
\end{align*}
is a (proper) submersion $W_\epsilon \to X \times
(a-\epsilon,a+\epsilon)$.  Thus the class $[W_\epsilon]$, as $\epsilon
\to 0$, is a well-defined element of $C_d^\pitchfork(X;a,a)$ and hence
gives an object
\begin{align*}
  \gamma(W,a) = ([W_\epsilon],a,a) \in \ob C_d^\pitchfork(X).
\end{align*}
This defines the functor $\gamma: D_d^\pitchfork \to C_d^\pitchfork$
on the level of objects, and it is defined similarly on morphisms.

\begin{proposition}\label{prop:5}
  The forgetful map $\alpha: \beta \ddp \to D_d$ is a weak
  equivalence.
\end{proposition}
\begin{proof}
  We apply the relative surjectivity criteria~\ref{cri:1} to the map
  $\beta \ddp \to D_d$. The argument is completely analogous to the
  proof of Proposition~4.2.4 of \cite{arXiv1}.
  
  First we show that $\beta\ddp(X) \to D_d(X)$ is surjective. Let $W
  \subseteq X\times \setR \times\setR^{d-1+\infty}$ be an element of
  $D_d(X)$.  For each $x \in X$ we can choose $a_x \in \setR$ such
  that $a_x$ is a regular value of $f_x: W_x = \pi^{-1}(x) \to \setR$.
  The same number $a_x$ will be a regular value of $f_y: W_y \to
  \setR$ for all $y$ in a small neighborhood $U_x \subseteq X$ of $x$.
  Therefore we can pick a locally finite open covering $\mathcal{U} =
  (U_j)_{j\in J}$ of $X$, and real numbers $a_j$, so that $f_j \colon
  W_j \to\setR$ is fiberwise transverse to $a_j$, where $W_j = W|U_j
  \in D_d(U_j)$.  Thus $(W_j, a_j)$ is an object of $\ddp(U_j)$ with
  $a_j\colon U_j\to\setR$ the constant map.

  For each finite subset $R\subseteq J$, set $W_R = W|U_R$ and $a_R =
  \min\{ a_j \mid j\in R\}$. If $R\subseteq S$ then $a_S \leq a_R$ and
  $(W_S, a_S, a_R)$ is an element $\phi_{RS} \in N_1\ddp(U_S)$. The
  pair $(\mathcal{U}, \Phi)$ with $\Phi = (\phi_{RS})_{R\subseteq S}$
  is an element of $\beta\ddp(U_S)$ that maps to $W$ by~$\alpha$.

  Second, let $A$ be a closed subset of $X$, $W \subseteq X \times
  \setR\times\setR^{d-1+\infty}$ an element of $D_d(X)$, and suppose
  we are given a lift to $\beta\ddp(U')$ of the restriction of $W$ to
  some open neighborhood $U'$ of $A$.  This lift is given by a locally
  finite open cover $\mathcal{U}' = \{ U_j | j \in J\}$, together with
  smooth functions $a_R: U_R \to \setR$, one for each finite non-empty
  $R \subseteq J$.  Let $J' \subseteq J$ denote the set of $j$ for
  which $U_j$ is non-empty, and let $J'' = J - J'$.

  Choose a smooth function $b: X \to [0,\infty)$ with $A\subseteq \Int
  b^{-1}(0)$ and $b^{-1}(0) \subseteq U'$.  Let $q = 1/b : X \to
  (0,\infty]$.  We can assume that $q(x) > a_R(x)$ for $R\subseteq J'$
  (make $U'$ smaller if not).  For each $x \in X - U'$, we can choose
  an $a \in \setR$ satisfying
  \begin{enumerate}[(i)]
  \item $a > q(x)$
  \item $a$ is a regular value for $f_x: \pi^{-1}(x) \to \setR$.
  \end{enumerate}
  The same number $a$ will satisfy (i) and (ii) for all $x$ in a small
  neighborhood $U_x \subseteq X - A$ of $X$, so we can pick an open
  covering $\mathcal{U}'' = \{U_j \mid j \in J''\}$ of $X - U'$, and
  real numbers $a_j$, such that (i) and (ii) are satisfied for all $x
  \in U_j$.  The covering $\mathcal{U}''$ can be assumed locally
  finite.  For each finite non-empty $R \subseteq J''$, set $a_R =
  \min\{ a_j \mid j \in R\}$.  For $R \subseteq J = J' \cup J''$,
  write $R = R' \cup R''$ with $R' \subseteq J'$ and $R'' \subseteq
  J''$, and define $a_R = a_{R'}$ if $R' \neq \emptyset$.

  This defines smooth functions $a_R: U_R \to \setR$ for all finite
  non-empty subsets $R\subseteq J$ ($a_R$ is a constant function for
  $R\subseteq J''$) with the property that $R\subseteq S$ implies $a_S
  \leq a_R|U_S$.  This defines an element of $\beta\ddp(X)$ which
  lifts $W \in D_d(X)$ and extends the lift given near $A$.
\end{proof}

\begin{proposition}\label{prop:B-gamma}
  The inclusion functor $\gamma: \ddp \to \cdp$ induces an equivalence
  $B|\ddp| \to B|\cdp|$.
\end{proposition}
\begin{proof}
  We show that $\gamma$ induces an equivalence $|N_k\ddp| \to
  |N_k\cdp|$ for all $k$, using the relative surjectivity
  criteria~\ref{cri:1}.

  An element of $N_k\cdp(X)$ can be represented by a sequence of
  functions $a_0 \leq \dots \leq a_k : X \to \setR$, a function
  $\epsilon: X \to (0,\infty)$, and a submanifold $W\subseteq X \times
  (a_0 - \epsilon, a_k + \epsilon) \times \setR^{d-1+\infty}$.
  Choosing a diffeomorphism $X \times (a_0 - \epsilon, a_k + \epsilon)
  \to X \times \setR$ which is the inclusion map on $X \times (a_0 -
  \epsilon/2, a_k + \epsilon/2)$, lifts the element to $N_k\ddp(X)$.
  This proves the absolute case and the relative case is similar.
\end{proof}

\begin{proposition}\label{prop:delta}
  The forgetful functor $\delta: C_d \to \cdp$ induces a weak equivalence
  $B|C_d| \to B|\cdp|$.
\end{proposition}
\begin{proof}
  Again we prove the stronger statement that $\delta$ induces an
  equivalence $|N_kC_d| \to |N_k\cdp|$ for all $k$.

  First, remember that two smooth maps $f\colon M\to P$ and $g\colon
  N\to P$ are called transversal if their product is transverse to the
  diagonal in $P\times P$.  We apply Criteria~\ref{cri:1}, and first
  prove that $\delta$ is surjective on concordance classes. Let
  $\psi\colon\setR\to [0, 1]$ be a fixed smooth function which is 0
  near $(-\infty,\frac13]$ and is 1 near $[\frac23,\infty)$,
  satisfying that $\psi' \geq 0$ and that $\psi' > 0$ on
  $\psi^{-1}((0,1))$.

  Given smooth functions $a_0 \leq a_1\colon X \to\setR$ with $(a_1 -
  a_0)^{-1}(0) \subseteq X$ an open subset, we define $\phi\colon
  X\times\setR \to X\times\setR$ by the formulas
  \begin{align*}
    \phi(x, u) &= (x, \phi_x(u)),\\
    \phi_x(u) &=
    \begin{cases}
      a_0(x) + (a_1(x) - a_0(x))
      \psi\bigl(\frac{u-a_0(x)}{a_1(x)-a_0(x)}\bigr) & \text{if $a_0(x)
        < a_1(x)$},\\
      a_0(x) & \text{if $a_0(x) = a_1(x)$}.
    \end{cases}
  \end{align*}

  Suppose that $W \in \cdp(X; a_0, a_1)$ with $a_0\leq a_1$.  The
  fiberwise transversality condition (iii) of Definition~\ref{def:3}
  implies that $(\pi, f)$ and $\phi$ are transverse, and hence that
  \begin{align*}
    W_\phi = \phi^* W = \{ (x, u, z) \mid \pi(z) = x, f(z) = \phi_x(u)
    \}
  \end{align*}
  is a submanifold of $X\times \setR\times W$. Using the embedding
  $W\subset X\times\setR\times\setR^{d-1+\infty}$ we can rewrite
  $W_\phi$ as
  \begin{align*}
    W_\phi &= \{ (x, u, r) \mid (x, \phi_x(u), r) \in W\} \subseteq
    X\times\setR\times\setR^{d-1+\infty}.
  \end{align*}
  It follows that
  \begin{align*}
    W_\phi \cap \bigl(X\times (-\infty, a_0 + \epsilon) \times
    \setR^{d-1+\infty}\bigr) &= M_0\times (-\infty, a_0+\epsilon)\\
    W_\phi \cap \bigl(X\times (a_1 - \epsilon, +\infty)
    \times\setR^{d-1+\infty}\bigr) &= M_1\times (a_1-\epsilon, +\infty),
  \end{align*}
  where $\epsilon = 1$ on $(a_1 - a_0)^{-1}(0)$ and $\epsilon =
  \frac13(a_1 - a_0)$ otherwise.  Thus $W_\phi$ defines an element of
  $C_d(X; a_0, a_1, \epsilon)$, and in turn an element of $C_d(X; a_0,
  a_1)$.

  We have left to check that $W_\phi$ is concordant to $W$ in $\cdp(X;
  a_0, a_1)$. To this end we interpolate between the identity and our
  fixed function $\psi\colon\setR\to [0, 1]$. Define
  \begin{align*}
    \psi_s(u) &= \rho(s) \psi(u) + (1-\rho(s))u
  \end{align*}
  with $\rho$ any smooth function from $\setR$ to $[0, 1]$ for which
  $\rho = 0$ near $(-\infty, 0]$ and $\rho=1$ near $[1, \infty)$.
  Define $\Phi\colon X\times\setR\times\setR\to X \times \setR$ as
  $\Phi(x,s,u) = (x,\Phi_x(s,u))$ where
  \begin{align*}
    \Phi_x(s, u) &=
    \begin{cases}
      a_0(x) + (a_1(x) - a_0(x))
      \psi_s\bigl(\frac{u-a_0(x)}{a_1(x)-a_0(x)}\bigr) & \text{if $a_0(x)
        < a_1(x)$},\\
      \rho(s)a_0(x) + (1-\rho(s)) u & \text{if $a_0(x) = a_1(x)$}.
    \end{cases}
  \end{align*}
  $\Phi$ is transversal to $(\pi,f)$, and the manifold
  \begin{align*}
    W_\Phi &= \{ ((x, s), u, r) \mid (x,\Phi_x(s, u), r) \in W\}
    \subseteq (X \times \setR) \times \setR \times \setR^{d-1+\infty}
  \end{align*}
  defines the required concordance in $\cdp(X\times\setR)$ from
  $W$ to $W_\phi$.

  We have proved that
  \begin{align*}
    \delta\colon N_0C_d[X] &\to N_0\cdp[X] \quad \text{and}\\
    \delta\colon N_1C_d[X] &\to N_1\cdp[X]
  \end{align*}
  are both surjective. The obvious relative argument is similar, and
  we can use Criteria~\ref{cri:1}. This proves that $\delta\colon
  \abs{N_k C_d} \to \abs{N_k \cdp}$ is a weak homotopy equivalence for
  $k = 0$ and $k = 1$. The case of general $k$ is similar.
\end{proof}
\begin{remark}
  \label{remark:a-locally-constant}
  There are versions of the sheaves $\ddp$, $\cdp$, $C_d$, where the
  functions $a: X \to \setR$ are required to be \emph{locally
    constant}.  The proofs given in this section remain valid for
  these sheaves (the point is that in the proof of
  Proposition~\ref{prop:5}, we are choosing the functions $a_j : U_j
  \to \setR$ locally constant anyway).  This proves the claim in the
  last sentence of Remark~\ref{rem:1}(ii).
\end{remark}


\section{Tangential structures}
\label{sec:tang-struct}

We prove the version of the Main Theorem with tangential structures,
as announced in the introduction.  First we give the precise
definitions.

Fix $d \geq 0$ as before, and let $BO(d) = G(d,\infty)$ denote the
Grassmannian of $d$-planes in $\setR^\infty$, $U_d \to BO(d)$ the
universal $d$-dimensional vector bundle, and $EO(d)$ its frame bundle.
Let
\begin{align*}
  \theta: B \to BO(d)
\end{align*}
be a Serre fibration (e.g.\ a fiber bundle).  We think of $\theta$ as
\emph{structures} on $d$-dimensional vector bundles: If $f: X \to
BO(d)$ classifies a vector bundle over $X$, then a $\theta$-structure
on the vector bundle is a map $l: X \to B$ with $\theta \circ l = f$.

An important class of examples comes from group representations.  If
$G$ is a topological group and $\rho: G \to \GL(d,\setR)$ is a
representation, then it induces a map $B\rho: BG \to B\GL(d,\setR)
\simeq BO(d)$, which we can replace by a Serre fibration.  In this
case, a $\theta$-structure is equivalent to a lifting of the structure
group to $G$.  These examples include $SO(d)$, $\mathrm{Spin}(d)$,
$\mathrm{Pin}(d)$, $U(d/2)$ etc.

Another important class of examples comes from spaces with an action
of $O(d)$.  If $Y$ is an $O(d)$-space, we let $B = EO(d) \times_{O(d)}
Y$.  If $Y$ is a space with trivial $O(d)$-action, then a
$\theta$-structure amounts to a map from $X$ to $Y$.  If $Y =
(O(d)/SO(d)) \times Z$, with trivial action on $Z$, then a
$\theta$-structure amounts to an orientation of the vector bundle
together with a map from $X$ to $Z$.

The proof of the main theorem applies almost verbatim if we add
$\theta$-structures to the tangent bundles of all $d$-manifolds in
sight.  We give the necessary definitions.

If $V\to X$ and $U\to Y$ are two vector bundles, a \emph{bundle map}
$V \to U$ is a continuous map of the total spaces of the vector
bundles, which on each fiber of $V$ restricts to a linear isomorphism
onto a fiber of $U$.  Let $\Bun(V,U)$ denote the space of all bundle
maps, equipped with the compact-open topology.  If $U = U_d$ is the
universal bundle over $BO(d)$
we have the following well known property.
\begin{lemma}
  Let $V \to X$ be a $d$-dimensional vector bundle with $X$
  paracompact.  Let $U_d \to BO(d)$ be the universal bundle.  Then the
  space $\Bun(V,U_d)$ is contractible.
\end{lemma}
\begin{proof}
  Since $U_d \subseteq BO(d) \times \setR^\infty$, we have a map
  \begin{align*}
    \Bun(V,U_d) \subseteq \Map(V,U_d) \to \Map(V,\setR^\infty)
  \end{align*}
  which identifies $\Bun(V,U_d)$ with the space of continuous maps $V
  \to \setR^\infty$ which restrict to linear monomorphisms on each
  fiber of $V \to X$.  Now define linear monomorphisms $\setR^\infty
  \to \setR^\infty$ by
  \begin{align*}
    i_1(x) &= (x_0, 0 , x_1, 0, x_2, 0, \dots)\\
    i_t(x) &= (1-t)x + ti_1(x), \quad \text{$0 \leq t\leq 1$}\\
    j(x) &= (0,x_0, 0, x_1, 0,x_2, \dots)
  \end{align*}
  There is an induced homotopy
  \begin{align}\label{eq:29}
    \begin{aligned}[]
      [0,1] \times \Map(V,\setR^\infty) &\to \Map(V, \setR^\infty)\\
      (t,f) &\mapsto i_t \circ f
    \end{aligned}
  \end{align}
  which restricts to a homotopy of self-maps of $\Bun(V,U_d)$,
  starting at the identity.

  It is well known that $\Bun(V,U_d)$ is non-empty when $X$ is
  paracompact (see e.g.\ \cite{MR0440554}).  Pick $g \in \Bun(V,U_d)$,
  and define a homotopy by
  \begin{align*}
    [0,1] \times \Map(V,\setR^\infty) &\to \Map(V,\setR^\infty)\\
    (t,f) & \mapsto (1-t) (i_1 \circ f) + t (j\circ g).
  \end{align*}
  This restricts to a homotopy of self-maps of $\Bun(V,U_d)$ which
  starts at $f \mapsto i_1 \circ f$.  Combined with the
  homotopy~\eqref{eq:29} we get a homotopy of self-maps of
  $\Bun(V,\setR^\infty)$ which starts at the identity and ends at the
  constant map to $j \circ g$.
\end{proof}

A (non-identity) point in $\mor\mathcal{C}_d$ is given by $(W, a_0,
a_1)$, where $a_0 < a_1\in \setR$ and $W$ is a submanifold (with
boundary) of $[a_0,a_1]\times \setR^{d-1+n}$, $n\gg 0$.  The tangent
spaces $T_pW$ define a map
\begin{align*}
  \tau_W: W \to G(d,n) \to BO(d),
\end{align*}
covered by a bundle map $TW \to U_d$.

\begin{definition}\label{defn:5.1}
\begin{sloppypar}
  Let $\mathcal{C}_\theta$ be the category with morphisms $(W, a_0,
  a_1, l)$, where $(W, a_0,a_1) \in \mor\mathcal{C}_d$ and $l: W \to
  B$ is a map satifying $\theta \circ l = \tau_W$.  We topologize
  $\mor\mathcal{C}_\theta$ as in~\eqref{eq:7}, but with $B_\infty(W)$
  replaced with $B_\infty^\theta(W) = \Emb^\theta(W,[0,1]\times
  \setR^{d-1+\infty})/\Diff(W)$, where $\Emb^\theta$ is defined by the
  pullback square
  \begin{align}
    \label{eq:8}
    \begin{aligned}
    \xymatrix{ {\Emb^\theta(W,[0,1]\times \setR^{d-1+\infty})}
      \ar[r]\ar[d] &
      {\Bun(TW,\theta^*U_d)} \ar[d]^\theta\\
      {\Emb(W, [0,1] \times \setR^{d-1+ \infty})} \ar[r]^-{\tau_W} &
      {\Bun(TW,U_d).}  }
    \end{aligned}
  \end{align}
  The objects of $\mathcal{C}_\theta$ are topologized similarly.
\end{sloppypar}
\end{definition}

The space $\Bun(TW,U_d)$ is contractible, so the inclusion of the
fiber product in the product
\begin{align*}
  \Emb^\theta(W,[0,1]\times \setR^{d-1+\infty}) \to \Emb(W, [0,1]
  \times \setR^{d-1+ \infty}) \times \Bun(TW,\theta^*U_d)
\end{align*}
is a homotopy equivalence.  Dividing out the action of $\Diff(W)$ we
get a homotopy equivalence
\begin{align*}
  B_\infty^\theta(W) \xrightarrow{\simeq} E \Diff(W)
  \times_{\Diff(W)} \Bun(TW,\theta^*U_d).
\end{align*}
Thus, up to homotopy,
\begin{align}\label{eq:3}
  \ob\mathcal{C}_\theta &\simeq \coprod_M E\Diff(M)
  \times_{\Diff(M)} \Bun(\setR \times TM,\theta^* U_d), \\
  \label{eq:4}
  \mor\mathcal{C}_\theta &\simeq \coprod_W E\Diff(W)
  \times_{\Diff(W)} \Bun(TW,\theta^*U_d),
\end{align}
where $M$ runs over closed $(d-1)$-manifolds, one in each
diffeomorphism class, and $W$ runs over compact $d$-dimensional
cobordisms, one in each diffeomorphism class.  As before, $\Diff(W)
\simeq \Diff(W,\{\partial_\mathrm{in}\}, \{\partial_\mathrm{out}\})$
denotes the topological group of diffeomorphisms that restrict to
diffeomorphisms of the incoming and outgoing boundaries separately
(or to product diffeomorphisms on a collar).

The left hand side of the homotopy equivalence~\eqref{eq:4} is the
space of all morphisms in $\mathcal{C}_\theta$.  The space of
morphisms between two fixed objects can be determined similarly.  We
first treat the case $\theta = \mathrm{id}$.  Let $c_0 = (M_0, a_0)$
and $c_1 = (M_1, a_1)$ be two objects of $\mathcal{C}_d$, given by
real numbers $a_0 < a_1$, closed manifolds $M_\nu \subseteq
\setR^{d-1+\infty}$.  Let $W$ be a compact manifold and $h_0:[0,1)
\times M_0 \to W$ and $h_1: (0,1] \times M_1 \to W$ be collars as
in~\eqref{eq:6}.  Let
\begin{align*}
  \Emb^\partial(W,[0,1] \times \setR^{d-1+\infty}) \subseteq
  \Emb(W,[0,1] \times \setR^{d-1+\infty}) 
\end{align*}
be the subspace consisting of embeddings $j$ which satisfy $j\circ
h_0(t,x) = (t,x)$ for $t$ sufficiently close to $0$ and $j \circ
h_1(t,x) = (t,x)$ for $t$ sufficiently close to 1.  Let
$\Diff(W;\partial W) \subseteq \Diff(W)$ be the subgroup consisting of
diffeomorphisms that restrict to the identity on a neighborhood of
$\partial W$.  This subgroup acts on $\Emb^\partial(W,[0,1] \times
\setR^{d-1+\infty})$ and we let $B_\infty^\partial(W)$ be the orbit
space
\begin{align*}
  B_\infty^\partial(W) = \Emb^\partial(W,[0,1] \times
  \setR^{d-1+\infty})/\Diff(W;\partial W).
\end{align*}
Then, up to homeomorphism, the space of morphisms is
\begin{align*}
  \mathcal{C}_d(c_0,c_1) \cong \coprod_W B_\infty^\partial(W),
\end{align*}
where the disjoint union is over cobordisms $W$ from $M_0$ to $M_1$,
one in each diffeomorphism class relative to $M_0$ and $M_1$.  Since
$\Emb^\partial(W,[0,1] \times \setR^{d-1+\infty})$ is contractible, we
get the homotopy equivalence
\begin{align*}
  \mathcal{C}_d(c_0,c_1) \simeq \coprod_W B\Diff(W;\partial W).
\end{align*}
The case of a general $\theta: B \to BO(d)$ is handled similarly.  If
$l_0: M_0 \to B$ and $l_1: M_1 \to B$ are two maps satisfying $\theta
\circ l_\nu = \tau_{\setR \times M_\nu}$ and $c_\nu = (M_\nu, a_\nu,
l_\nu)$, $\nu=0,1$, then we get
\begin{align}\label{eq:9}
  \mathcal{C}_\theta(c_0,c_1) \simeq \coprod_W E\Diff(W;\partial W)
  \times_{\Diff(W;\partial W)} \Bun^\partial(TW,\theta^* U_d),
\end{align}
where $\Bun^\partial(TW,\theta^* U_d) \subseteq \Bun(TW, \theta^*
U_d)$ is the subspace consisting of bundle maps which agree with the
maps induced by $l_0$ and $l_1$ over a neighborhood of $\partial W$.

Let us consider the case of ordinary orientations in more detail.
Here $B=BSO(d)$ is the oriented Grassmanian consisting of
$d$-dimentional linear subspaces of $\setR^\infty$ together with a
choice of orientation, and $\theta:B\to BO(d)$ is the twofold covering
space that forgets the orientation.  Let $W$ be a cobordism between
the oriented manifolds $M_0$ and $M_1$.  Then the set
\begin{align*}
  \mathrm{Or}(W;\partial W) = \pi_0\Bun^\partial(TW,\theta^* U_d)
\end{align*}
is the set of orientations of $W$ agreeing with the orientations given
near $\partial W$ (i.e.\ the collars $h_0$ and $h_1$ are oriented
embeddings).  Furthermore, the connected components of
$\Bun^\partial(TW,\theta^* U_d)$ are contractible, so we get a
homotopy equivalence
\begin{multline*}
  E\Diff(W;\partial W) \times_{\Diff(W;\partial W)}
  \Bun^\partial(TW,\theta^* U_d) \simeq\\ E\Diff(W;\partial W)
  \times_{\Diff(W;\partial W)} \mathrm{Or}(W;\partial W).
\end{multline*}
The stabilizer of an element of $\mathrm{Or}(W;\partial W)$ is the
subgroup $\Diff^{+}(W;\partial W)$ of orientation preserving
diffeomorphisms, restricting to the identity near the boundary.  Thus
we get
\begin{align*}
  \mathcal{C}_d^+(c_0,c_1) \simeq \coprod_W B\Diff^{+}(W; \partial W),
\end{align*}
where the disjoint union is over all oriented cobordisms $W$ from
$M_0$ to $M_1$, one in each oriented diffeomorphism class.

\begin{definition}
  Let $\theta_{d,n}: B_{d,n} \to G(d,n)$ be the pullback
  \begin{align*}
    \xymatrix{
      B_{d,n} \ar[r]\ar[d]_{\theta_{d,n}} & B \ar[d]^\theta \\
      G(d,n) \ar[r] & BO(d),
    }
  \end{align*}
  and let $\G{\theta}$ be the spectrum whose $(n+d)$th space is
  $\mathrm{Th}(\theta_{d,n}^* U_{d,n})$.
\end{definition}

The cofiber sequence~\eqref{eq:54} generalizes to a cofiber sequence
\begin{align*}
  \G{\theta} \longrightarrow \Sigma^\infty B_+ \longrightarrow
  \G{\theta_{d-1}},
\end{align*}
where $\theta_{d-1}$ is the pullback
\begin{align*}
  \xymatrix{
    B_{d-1} \ar[r]\ar[d]_{\theta_{d-1}} & B \ar[d]^\theta \\
    BO(d-1) \ar[r] & BO(d).  }
\end{align*}

With these definitions, the general form of the main theorem (as also
stated in the introduction) is that for every tangential structure
$\theta$, there is a weak equivalence
\begin{align*}
  B \mathcal{C}_\theta \simeq \Omega^{\infty-1} \G{\theta} =
  \colim_{n \to \infty} \Omega^{d+n-1}
  \Th(\theta_{d,n}^*U_{d,n}^\perp).
\end{align*}

The $\theta$-versions of the sheaves used in section 4 to prove the
special case $\theta = \id$, are defined as follows.

\begin{definition}
  Let $W \in D_d(X)$.  Let $T^\pi W$ be the fiberwise tangent bundle
  of the submersion $\pi: W \to X$.  The embedding $W\subset X \times
  \setR^{d+\infty}$ induces a canonical classifying map $T^\pi W: W
  \to BO(d)$.  Let $D_\theta(X)$ be the set of pairs $(W,l)$ with $W
  \in D_d(X)$ and $l: W \to B$ a map satisfying $\theta \circ f =
  T^\pi W$.
\end{definition}

The sheaves $C_d$, $C_d^\pitchfork$ and $D_d^\pitchfork$ all consist
of submanifolds $W \subseteq X \times \setR^{d+n}$ such that the
projection $\pi: W \to X$ is a submersion, together with some extra
data.  The tangential structure versions $C_\theta$,
$C_\theta^{\pitchfork}$ and $D_\theta^{\pitchfork}$ are defined in the
obvious way: add a lifting $l: W \to B$ of the vertical tangent bundle
$T^\pi W: W \to BO(d)$.

With these definitions, the proofs of section
\ref{sec:proof-of-main-theorem} apply almost verbatim.  We note that
the $\theta$-versions of Theorem~\ref{thm:4} and
Proposition~\ref{prop:delta} use that $\theta$ is a Serre fibration.


\section{Connectedness issues}
\label{cha:connectedness-issues}

This section, technically the hardest of the paper, compares the
category $\mathcal{C}_\theta$ with the \emph{positive boundary}
subcategory $\mathcal{C}_{\theta,\partial}$.  It is similar in spirit
to section \S6 of \cite{arXiv1}.  The two categories have the same
space of objects. The space of morphisms of
$\mathcal{C}_{\theta,\partial}$ is as in~\eqref{eq:7} and
Definition~\ref{defn:5.1}, but taking only disjoint union over the $W$
for which each connected component has non-empty \emph{outgoing}
boundary: if $W$ is a cobordism from $M_0$ to $M_1$, then $\pi_0 M_1
\to \pi_0 W$ is surjective.  In this section we prove
\begin{theorem}\label{thm:connectedness}
  For $d \geq 2$ and any $\theta: B \to BO(d)$, the inclusion
  \begin{align*}
    B\mathcal{C}_{\theta,\partial} \to B\mathcal{C}_{\theta}
  \end{align*}
  is a weak equivalence.
\end{theorem}

In order to simplify the exposition we treat only the case $\theta =
\mathrm{id}$.  The general case of an arbitrary $\theta$-structure is
similar.

We say that a map $f: X \to Y$ of topological spaces is
\emph{$\pi_0$-surjective} if the induced map $\pi_0X \to \pi_0 Y$ is
surjective.  The subsheaf $\ddpb\subseteq \ddp$ is defined as follows:
$(W, a_0, a_1) \in \ddp(\point)$ is in $\ddpb(\point)$ if the
inclusion
\begin{align*}
  f^{-1} (a_1) \to f^{-1} [a_0, a_1]
\end{align*}
is $\pi_0$-surjective.  In general $\chi = (W, a_0, a_1) \in \ddp(X)$
is in $\ddpb(X)$ if $\chi_{|\{x\}} \in \ddpb(\{x\})$ for all $x \in
X$.  The proof given above that $|\beta D^\pitchfork_d| \simeq
B\mathcal{C}_d$ (in Propositions \ref{prop:B-eta}, \ref{prop:B-gamma}
and \ref{prop:delta}) is easily modified to show that $|\beta
D^\pitchfork_{d,\partial}| \simeq B\mathcal{C}_{d,\partial}$.  We will
show that the composite map of sheaves $\beta \ddpb \to \beta \ddp \to
D_d$ satisfies the relative lifting criteria~\ref{cri:1} for all
$d\geq 2$.

\subsection{Discussion}
\label{sec:discussion}

We describe the ideas involved and indicate the issues in proving
that the map $\beta D_{d,\partial}^\pitchfork \to D_d$ is a weak
equivalence.

As a first approximation we can try to repeat the proof for $\beta
D_d^\pitchfork \to D_d$ (in Proposition~\ref{prop:5}), by choosing
regular values $a_x \in \setR$ for $f_x: W_x \to \setR$ ``at random''
(using Sard's theorem), and using that $a_x$ is a regular value for
$f_y: W_y \to \setR$ also for $y$ in a small neighborhood $U_x$ of $x
\in X$.  This will produce an element $(W,(U_j,a_j)_{j \in J}) \in
\beta D_d^\pitchfork(X)$ but in general there is, of course, no reason
to expect to get an element of $\beta D_{d,\partial}^\pitchfork(X)
\subseteq \beta D_d^\pitchfork(X)$.  The idea is now to deform (i.e.\
change by a concordance) the underlying $W \in D_d(X)$ to an element
$W' \in D_d(X)$ such that $W'$ together with the regular values $a_j$
(possibly slightly perturbed) defines an element of $\beta
D_{d,\partial}^\pitchfork(X)$.

\begin{figure}[htbp]
  \centering
  \includegraphics[scale=0.8]{picture/GMTW.20}
  \caption{ }
  \label{fig:1a}
\end{figure}

It is instructive to first consider the case $X = \point$.  Given an
element $(W,a_0 < \dots < a_k) \in N_k D_d^\pitchfork(\point)$, it is
easy to see that there is a concordance $H \in D_d(\setR)$ from $W$ to
$W'$ such that $(W', a_0 < \dots < a_k) \in N_k
D_{d,\partial}^\pitchfork(\point)$.  Roughly, we have to get rid of
some \emph{local maxima}, with values between $a_0$ and $a_k$, of the
function $f: W \to \setR$ cf.\ Figure~\ref{fig:1a}.  
 A naive way to do that is to ``pull them up'', i.e.\
if $p \in W$ is near a ``local maximum'' for $f: W \to \setR$, then we
can change $f$ near $p$ to have $f(p) > a_k$ cf.\ Figure~\ref{fig:1b}.
A better way (for reasons explained below) to get rid of a local
maximum, is given in Lemma~\ref{lemma:surgery} below.

\begin{figure}[htbp]
  \centering
  \includegraphics[scale=0.8]{picture/GMTW.21}
  \caption{ }
  \label{fig:1b}
\end{figure}

For general $X$ it is equally easy to solve the problem
\emph{locally}.  Given $W \in D_d(X)$, suppose we have chosen regular
values $a_j \in \setR$ and corresponding open covering $U_j \subseteq
X$, $j \in J$, such that $(W, (a_j,U_j)_{j \in J})$ defines an element
of $\beta D_d^\pitchfork(X)$.  Given $x \in X$ it is easy (as in the
case $X = \point$) to find a small neighborhood $U_x \subseteq X$ and
a concordance $H_x \in D_d(U_x \times \setR)$ from $W | U_x$ to $W'
\in D_d(U_x)$ such that $(W', (a_j,U_j\cap U_x)_{j \in J})$ defines an
element of $\beta D_{d,\partial}^\pitchfork(U_x)$.  We now need to
\emph{glue} these local constructions.

The locally defined concordance $H_x \in D_d(U_x \times \setR)$ can be
assumed to extend to $H_x \in D_d(X \times \setR)$.  Namely we may
choose a bump function $\lambda: X \to [0,1]$, supported in $U_x$, and
which is 1 in a smaller neighborhood $U_x' \subseteq U_x$, and let $h:
U_x \times \setR \to U_x \times \setR$ be given by $h(x,t) =
(x,t\lambda(x))$.  Then $H_x' = h^* H_x \in D_d(U_x \times \setR)$ is
a concordance which is constant outside the support of $\lambda$, so
it extends to a concordance $H_x' \in D_d(X \times \setR)$.  Moreover
$H_x' | (U_x' \times \setR) = H_x | (U_x' \times \setR)$.  Thus $H_x'$
is a concordance from $W$ to $W' \in D_d(X)$, such that $W'|U_x' \in
D_d(U_x')$ lifts to $\beta D_{d,\partial}(U_x')$.  Also $W$ and $W'$ agree
outside $U_x \supseteq U_x'$.

We have described how, given a way of getting rid of a single local
maxima, to deform an element $W \in D_d(X)$
into an element $W'\in D_d(X)$, with the property that $W'|U_x'$ lifts
to $\beta D_{d,\partial}^\pitchfork(U_x')$, and such that $W$ and $W'$ agree
outside a larger open neighborhood $U_x \supseteq U_x'$.  Roughly, the
idea is now to apply such a construction for sufficiently many $x \in
X$, enough that the sets $U_x'$ cover $X$.  For this to work there is
one critical issue, however.  Namely it is essential that the local
construction used to get rid of fiberwise local maxima over $U_x'$
does not create new fiberwise local maxima over $U_x - U_x'$.  Without
this, the idea to ``apply such a construction for sufficiently many $x
\in X$'' will not work.

The naive idea of ``pulling local maxima up'' will not work, precisely
for this reason.  If we ``pull up'' a fiberwise local maximum over
$U_x'$, we have to pull less and less over $U_x - U_x'$ (as specified
by the bump function $\lambda$), which will give rise to fiberwise
local maxima of $f': W' \to \setR$ over $U_x - U_x'$ which are not
fiberwise local maxima of $f: W \to \setR$.

Thus we will need a way of deforming $f: W \to \setR$ to get rid of
local maxima without creating new ones in the process.  Such a
construction is described in Lemma~\ref{lemma:surgery} below.  It
describes a family of maps $f_t: K_t \to \setR$, $t \in [0,1]$ from
$d$-manifolds $K_t$, such that $f_0$ is the constant map $0: \setR^d
\to \setR$, such that $f_1: \setR^d - \{0\} \to \setR$ has $\lim_{x
  \to 0} f(x) = + \infty$, and such that $f_t: K_t \to \setR$ has no
local maxima, except some with value $0 \in \setR$, for any $t \in
[0,1]$.  Moreover each $K_t$ contains the open subset $\setR^d - D^d
\subseteq K_t$ and $f_t | (\setR^d - D^d) = 0$.

\subsection{Surgery}
\label{sec:surgery}

The geometric construction is based on the following lemma.  Let us
say that a map $f: M \to N$ is proper \emph{relative to} an open set
$U \subseteq M$, if $f_{|M-U}: M-U \to N$ is proper.
\begin{lemma}\label{lemma:surgery}
  There exists a smooth $(d+1)$-manifold $K$ containing $U = \setR
  \times (\setR^d - D^d)$ as an open subset, and smooth maps $(\pi,f):
  K \to \setR \times \setR$, such that
  \begin{enumerate}[(i)]
  \item\label{item:5} $\pi$ is a submersion, and $(\pi,f)$ is proper
    relative to $U$.  In particular, if we let $K_t = \pi^{-1}(t)$ and
    $U_t = U \cap K_t = \{t\} \times (\setR^d - D^d)$, then $f_t: K_t
    \to \setR$ is proper relative to $U_t$.
  \item\label{item:4} $(\pi,f)(t,x) = (t,0)$ for all $(t,x) \in U
    \subseteq K$.
  \item\label{item:6} $K_0 = \{0\}\times \setR^d$ and $f_0: K_0 \to
    \setR$ is the zero function.
  \item\label{item:7} For all $t \in [0,1]$ and all $a_0 < a_1 \in
    \setR$, the following inclusions are $\pi_0$-surjections
    \begin{align*}
      U_t \amalg f_t^{-1}(a_1) & \to f_t^{-1} ([a_0,a_1])\quad
      \text{if $0 \in [a_0,a_1]$}\\
      f_t^{-1}(a_1) & \to f_t^{-1} ([a_0,a_1])\quad \text{if $0
        \not\in [a_0,a_1]$.}
    \end{align*}
  \item \label{item:9}For all $a_0 < a_1 \in \setR$, the inclusion
    \begin{align*}
      f_1^{-1}(a_1) \to f_1^{-1} ([a_0,a_1])
    \end{align*}
    is a $\pi_0$-surjection.
  \item \label{item:8} $K_1 = \{1\} \times (\setR^d - \{0\})$ and
    $f_1: K_1 \to \setR$ is non-negative and has $0 \in \setR$ as only
    critical value.
  \item \label{item:16} $T^\pi K$ is a trivial vector bundle.
  \end{enumerate}
\end{lemma}
The last property, that $T^\pi K$ be a trivial vector bundle, is needed
to make the constructions work in the presence of $\theta$-structures.

As stated, the lemma is true also for $d=1$, but is useful only for $d
> 1$.  For $d>1$ the set $U_t$ is connected, and the properties
(\ref{item:6}) and (\ref{item:7}) say that the number of elements in
the quotient
\begin{align*}
  Q_t = \pi_0(f_t^{-1}[a_0,a_1])/\pi_0(f_t^{-1}(a_1))
\end{align*}
is never larger than the number of elements in $Q_0$.  For $0 \in
[a_0,a_1]$ and $d>1$, the inclusion $U_t \to f_t^{-1}([a_0,a_1])$
defines an element $[U_t] \in Q_t$, and (\ref{item:9}) says that
$[U_0] \in Q_0$ is not the basepoint, then $Q_1$ is strictly smaller
than $Q_0$.

\begin{proof}
  We will construct $K$ as a certain pullback of a 2-manifold $L$
  which we first construct.  $L$ will come with an immersion $(\pi,j):
  L \to [0,4] \times [0,\infty)$ and a function $f: L \to \setR$.  $L$
  will be glued from four pieces $L^1, \dots, L^4$ which we construct
  individually.  The pieces $L^1$, $L^2$ and $L^4$ will be subsets of
  $[0,1] \times [0,\infty)$, and $L^3$ will be the disjoint union of
  three open subsets of $[0,1] \times [0,\infty)$.  In all cases,
  $(\pi,j): L^\nu \to [0,1] \times [0,\infty)$ will be given by the
  inclusions.

  Let $\rho: [0,\infty) \to [0,1]$ be a smooth function with
  $\supp(\rho) = [0,1]$, $\rho(0) = 1$, and $\rho'(r) \leq 0$.  For $s
  \in [0,1]$ let $q_s(r) = \rho(r^2) \frac{1-s}{r^2 + s}$ and let
  $g_s$ and $\hat g_s$ be the functions given by
  \begin{align*}
    g_s(r) &= -q_s(r) - q_s(r-2) + q_0(r-1)\\
    \hat g_s (r) &= \sgn(r(r-2))\!\left( -q_0(r) - q_0(r-2) +
    q_{1-s}(r-1) - \frac{1-s}s \right) + \frac{1-s}s.
  \end{align*}
  $g_s(r)$ is defined unless $r = 1$ or $(s,r) \in
  \{0\}\times\{0,2\}$.  $\hat g_s(r)$ is defined unless $r \in
  \{0,2\}$ or $(s,r) = (1,1)$ or $(s,r) \in \{0\} \times [0,2]$.  It
  is easily checked that $g_s'(r) = 0$ only if $r\geq 3$, if $(s,r)
  \in (0,1] \times \{0,2\}$, or if $(s,r) \in \{1\} \times
  [2,\infty)$.  Similarly $\hat g_s'(r) = 0$ only if $r\geq 3$ or
  $(s,r) \in (0,1)\times \{1\}$.  All isolated critical points of
  $g_s$ and $\hat g_s$ are local minima.

  Define functions $f^\nu: L^\nu \to \setR$ for $\nu = 1,2,4$
  by the following formulas, using the (calculus) convention that the
  set $L^\nu \subseteq [0,1] \times [0,\infty)$ is the largest open
  set for which the definitions make sense.
  \begin{align*}
    f^1(t,r) & = \hat g_0(r + 3(1-t))\\
    f^2(t,r) & = \hat g_t(r)\\
    f^4(t,r) & = g_t(r+t).
  \end{align*}
  To define $f^3$, let $L^3 = L^3_- \amalg L^3_+ \amalg L^3_0$, where
  \begin{align*}
    L^3_- &= \{(t,r)\in [0,1]\times [0,\infty) \mid t < r < t+1 \}\\
    L^3_+ &=  \{(t,r)\in [0,1]\times [0,\infty) \mid (1-t) < r < (2-t) \} \\
    L^3_0 &= [0,1] \times (2,\infty).
  \end{align*}
  Let $f^3 = f^3_- \amalg f^3_+ \amalg f^3_0$, where
  \begin{align*}
    f^3_\epsilon(t,r) = \hat g_1(r + \epsilon t).
  \end{align*}

  It is easily checked that $f^1(1,r) = f^2(0,r)$, $f^2(1,r) =
  f^3(0,r)$ and $f^3(1,r) = f^4(0,r)$, so they glue to a continuous
  function $\tilde f: \tilde L \to \setR$, where $\tilde L$ is glued
  from $L^1, \dots, L^4$.  $\tilde L$ is a smooth manifold and comes
  with an immersion $(\tilde \pi, \tilde j): \tilde L \to [0,4] \times
  [0,\infty)$.  The 2-manifold $\tilde L$ is sketched in
  Figure~\ref{fig:1}, 
  which also depicts the map $\tilde \pi: \tilde L \to
  [0,4]$ as the projection onto the horizontal axis and $\tilde j: \tilde L
  \to [0,\infty)$ as the projection onto the vertical axis.

  The function $\tilde f$ is not smooth in the $t$-variable along the
  gluing lines.  To fix that, we choose a function $\sigma: [0,4] \to
  [0,4]$ which for each $n = 1,2,3$ has $\sigma(t) = n$ for all $t$
  near $n$.  Then let $L$ be defined by the pullback diagram
  \begin{align*}
    \xymatrix{
      L \ar[r]^{\overline \sigma} \ar[d]_\pi & {\tilde L} \ar[d]^{\tilde \pi}\\
      [0,4] \ar[r]^\sigma & [0,4],  }
  \end{align*}
  and let $j = \tilde j \circ \overline \sigma: L \to [0,\infty)$ and
  $f = \tilde f \circ \overline \sigma: L \to \setR$.  The resulting
  $f: L \to \setR$ is then smooth.

  Let $\lambda: \setR \to [0,1]$ be a smooth function which is $0$
  near $(-\infty,0]$ and 1 near $[1,\infty)$ and has $\lambda' > 0$ on
  $\lambda^{-1}((0,1))$.  Let $g: \setR\times \setR^d \to [0,4] \times
  [0,\infty)$ be the map given by
  \begin{align*}
    g(t,x) = (4 \lambda(t), 3 \lvert x \rvert^2).
  \end{align*}
  \begin{figure}[htbp]
    \centering
    \includegraphics[scale=0.7]{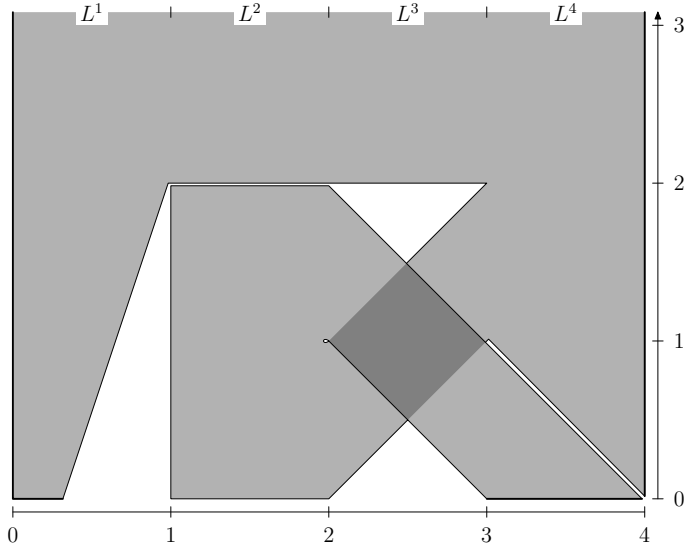}
    \caption{Image of $(\tilde \pi,\tilde j): \tilde L \to
      [0,4]\times[0,\infty)$.}
    \label{fig:1}
  \end{figure}

  \noindent
  To construct the map $(\pi,f): K \to \setR \times \setR$ of the
  proposition, define $K$ as the pullback in the diagram
  \begin{align}
    \label{eq:diagram6.1}
    \vcenter{
    \xymatrix{
      K \ar[d]_{(\pi,j)} \ar[r] & L \ar[d]^{(\pi,j)} \ar[r]^f &
      {\setR}\\
      {\setR \times \setR^d} \ar[r]^-g & {[0,4]\times [0,\infty).} &
    }}
  \end{align}
  Then $(\pi,j): K \to \setR \times \setR^d$ is a codimension 0
  immersion, and over $U = \setR \times (\setR^d - D^d)$ it is a
  diffeomorphism.  The diagram also provides a map $f: K \to \setR$,
  and it is easily seen that $(\pi,f): K \to \setR \times \setR$
  satisfies the first six properties of the proposition.  The
  differential of $(\pi,j): K \to \setR \times \setR^d$ defines a
  trivialization of the $d$-dimensional vector bundle $T^\pi K$.
\end{proof}

\begin{figure}[htbp]
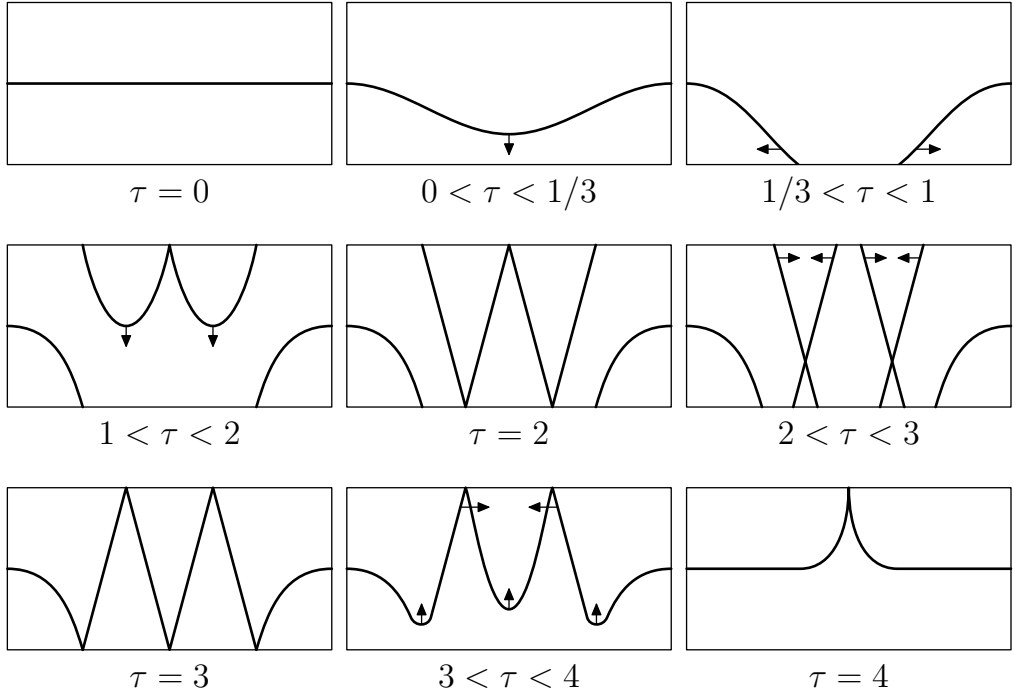

  \centering
  \includegraphics[width=.3\textwidth]{picture/GMTW.11}
  \includegraphics[width=.3\textwidth]{picture/GMTW.12}
  \includegraphics[width=.3\textwidth]{picture/GMTW.13}
  \\[3mm]
  \includegraphics[width=.3\textwidth]{picture/GMTW.14}
  \includegraphics[width=.3\textwidth]{picture/GMTW.15}
  \includegraphics[width=.3\textwidth]{picture/GMTW.16}
  \\[3mm]
  \includegraphics[width=.3\textwidth]{picture/GMTW.17}
  \includegraphics[width=.3\textwidth]{picture/GMTW.18}
  \includegraphics[width=.3\textwidth]{picture/GMTW.19}
  \caption{$(f_t,j_t)(K_t)$ for $d=1$ and various values of $\tau =
    \sigma(4\lambda(t))\in[0,4]$.
      }
  \label{fig:2}
\end{figure}

The manifold $K$ and the map $(\pi,f): K \to
\setR^d \times \setR$ are illustrated in Figure~\ref{fig:2}, which
shows the $d$-manifold $K_t = \pi^{-1}(t)$ for $d=1$ and various
values of $t \in[0,1]$.  The horizontal axis is $[-1,1] = D^d
\subseteq \setR^d$ and the projection is the immersion $j_t: K_t \to
\setR^d$.  The vertical axis is $(-\infty,\infty)$ and the projection
is the function $f_t: K_t \to (-\infty,\infty)$.  The small arrows
indicate how $K_t$ changes when $t$ increases.

Given an element $W \in D_d(\point)$, assume $e: \setR^d \to W$ is an
embedding with $e(\setR^d) \subseteq f^{-1}(r)$ for some $r \in
\setR$.  Then $W \times \setR \in D_d(\setR)$ has an embedded $\setR^d
\times \setR$ from which we can remove $D^d \times \setR$ and glue in
the manifold $K$ from the above Lemma~\ref{lemma:surgery} along the
embedded $(\setR^d - D^d) \times \setR$.  This gluing is over $\setR$
if we equip $K$ with the map $f+r: K \to \setR$ and we get a
concordance $W^e \in D_d(\setR)$ starting at $W \in D_d(\{0\})$.  We
will describe an enhanced version of this construction where we start
with $W\in D_d(X)$ and a finite set of embeddings
$e_\tau:X\times\setR^d\to W$ ($\tau\in T$) such that $r_\tau(x)=f\circ
e_\tau (x,u)$ is independent of $u\in\setR^d$.  The enhanced
construction will give an element $W^e\in D_d(X\times \setR^T)$ which
upon restriction to $X\times \{1\}^T$ is an element where the ``local
maxima'' at $e_\tau(x,0)$ have disappeared.

\begin{definition}\label{definition:E-r}
  Let $X$ be a manifold and $T$ a finite set.  Let $r: X \times T \to
  \setR$ be smooth.  For $\tau \in T$, let $q_{\tau,r}: (X \times
  \setR^T)\times \setR \to \setR \times \setR$ be the map
  \begin{align*}
    q_{\tau,r}((x,l),t) = (l_\tau, t - r(x,\tau)), \quad
    l=(l_\tau)_{\tau\in T}\ .
  \end{align*}
  Considering $K$ as a space over $\setR \times \setR$ via the map
  $(\pi,f)$ from \eqref{eq:diagram6.1}, we get a manifold
  $q_{\tau,r}^*K$ over $(X \times \setR^T) \times \setR$, containing
  $q_{\tau,r}^*U = (X \times \setR^T)\times (\setR^d - D^d)$ as an
  open subset.  Let
  \begin{align*}
    K^r &= \coprod_{\tau \in T} q_{\tau,r}^*K, & U^r &= \coprod_{\tau
      \in T} q_{\tau,r}^*U \subseteq K^r.
  \end{align*}
  This comes equipped with a map $(\pi^r, f^r): K^r \to (X \times
  \setR^T) \times \setR$ which is proper relative to $U^r = (X \times
  \setR^T) \times \coprod_T (\setR^d - D^d)$, and $\pi^r: K^r \to X
  \times \setR^T$ is a submersion.
\end{definition}
\begin{remark}
  \label{remark:E-monoidal}
  This behaves well under union in the $T$-variable.  If $T = T_0
  \amalg T_1$ and $r_\nu: X \times T_\nu \to \setR$, $\nu= 0,1$ are
  the restrictions of $r$, then
  \begin{align*}
    K^r = \proj_{X \times \setR^{T_0}}^*(K^{r_1}) \amalg \proj_{X
      \times \setR^{T_1}}^*(K^{r_0})
  \end{align*}
  where the indicated projections are $X \times \setR^T \to X \times
  \setR^{T_\nu}$, $\nu=0,1$.
\end{remark}

\begin{construction}\label{constr:1}
  Let $W \in D_d(X)$, and let $T$ be a finite set.  Let $r: X \times T
  \to \setR$ be smooth.  Then $X \times \coprod_T \setR^d = X \times T
  \times \setR^d$ is a space over $X \times \setR$ via the projection
  composed with $r$.  Let
  \begin{align*}
    e: X \times \coprod_T \setR^d \to W
  \end{align*}
  be an embedding over $X \times \setR$, i.e.\ with $\pi\circ
  e(x,\tau,u)=x$ and $f\circ e(x,\tau,u)=r(x,\tau)$.  This induces an
  embedding
  \[
  \tilde e: (X
  \times \setR^T) \times \coprod_T \setR^d \to \proj_X^*W,
  \]
  where
  $\proj_X: X \times \setR^T \to X$ is the projection.  Let $W^e$ be
  the pushout
  \begin{align}\label{eq:28}
    \begin{aligned}
    \xymatrix@C=1.5cm{
      U^r\ar[r]^-{\tilde e} \ar[d] & {\proj_X^* W - \tilde e(X \times
        \setR^T \times \coprod_T D^d)} \ar[d]\\
      K^r \ar[r] & W^e
    }
    \end{aligned}
  \end{align}
  This gives a manifold $W^e$ over $(X \times \setR^T) \times \setR$
  which defines an element of $D_d(X \times \setR^T)$.
\end{construction}

Elements of $D_d(X \times \setR^T)$ are submanifolds of $(X \times
\setR^T) \times \setR \times \setR^{d-1+\infty}$, so strictly speaking
the construction of $W^e$ includes a choice of an embedding
\begin{align*}
  \phi: W^e \to (X \times \setR^T) \times \setR \times
  \setR^{d-1+\infty}
\end{align*}
extending the given map $W^e \to (X \times \setR^T) \times \setR$.
Then the image $\phi(W^e)$ is an element of $D_d(X \times \setR^T)$.
The element $\proj_X^*W \in D_d(X \times \setR^T)$ has a preferred
embedding $i: \proj_X^*W \to (X \times \setR^T) \times \setR \times
\setR^{d-1+\infty}$ (namely the inclusion), and it is convenient to
assume that $\phi$ and $i$ agree on the subspace ${\proj_X^* W -
  \tilde e(X \times \setR^T \times \coprod_T D^d)}$.  Such an
embedding $\phi$ can always be chosen, and is unique up to isotopy.
It is irrelevant for the arguments which $\phi$ we choose, and
therefore we omit it from the notation, writing $W^e \in D_d(X \times
\setR^T)$ instead of $\phi(W^e)$.

\subsection{Connectivity}

We will apply the surgery construction of the previous section to a
morphism $(W,a_0,a_1) \in \ddp(X)$ with $a_0<a_1$.  The resulting $W^e
\in D_d(X \times \setR^T)$ will usually \emph{not} give rise to an
element $(W^e, a_0, a_1) \in \ddp(X \times \setR^T)$ because $f^e: W^e
\to \setR$ might not be fiberwise transverse to $a_0, a_1$.  Let $V =
V(a_0, a_1) \subseteq X \times \setR^T$ be the open set of points
$(x,l)$ for which $f^e_{(x,l)}: W^e_{(x,l)} \to \setR$ is transverse
to $a_0(x)$ and $a_1(x)$.  Then we have $(W^e, a_0, a_1)_{|V} \in
\ddp(V)$.  By Sard's theorem, any $(x,t) \in X \times \setR^T$ is in
$V(b_0,b_1)$ for some $b_0, b_1$ arbitrarily close to $a_0, a_1$.  The
goal is to use these concordances to get an element of $\ddpb$.  Since
the condition for being in $\ddpb\subseteq \ddp$ is pointwise, we
restrict attention to the case $X = \point$ in the following
propositions.
\begin{proposition}\label{proposition:connect}
  Let $(W, a_0, a_1) \in \ddp(\point)$ with $a_0<a_1$.  Let $r: T \to
  \setR$ and $e: \coprod_T \setR^d \to W$ be as in
  Construction~\ref{constr:1}.  Let $V = V(a_0, a_1)\subseteq \setR^T$
  be as above.
  \begin{enumerate}[(i)]
  \item\label{item:3} If $r(\tau) \neq a_0, a_1$ for all $\tau \in T$, then
    $\{0,1\}^T \subseteq V$.
  \item\label{item:15} If $(W, a_0, a_1) \in \ddpb(\point)$, then $(W^e, a_0,
    a_1)_{|V} \in \ddpb(V)$.
  \item\label{item:14} If $(W,a_0, a_1) \in \ddp(\point)$, $a_0< r<
    a_1$, and if
    \begin{align*}
      f^{-1}(a_1) \amalg \coprod_T \setR^d \to f^{-1}([a_0,a_1])
    \end{align*}
    is $\pi_0$-surjective, then the restriction to $\{1\}^T \subseteq
    \setR^T$ defines an element $(W, a_0, a_1)_{\{1\}^T} \in
    \ddpb(\{1\}^T)$.
  \end{enumerate}
\end{proposition}
\begin{proof}
  Let $l \in \{0,1\}^T$.  By Lemma~\ref{lemma:surgery}(\ref{item:8})
  we get that critical values of $f_l: W^e_l \to \setR$ will be either
  critical values of $f: W \to \setR$, or values $r(\tau)$ for
  $\tau\in T$ with $l_\tau = 1$.  This proves \emph{(\ref{item:3})}.
  \emph{(\ref{item:15})} follows from
  Lemma~\ref{lemma:surgery}\emph{(\ref{item:7})} and
  \emph{(\ref{item:14})} follows in the same way from
  Lemma~\ref{lemma:surgery}\emph{(\ref{item:9})}.
\end{proof}

If not $l \in V(a_0,a_1)$, then $l \in V(b_0, b_1)$ for some $b_0,
b_1$ near $a_0, a_1$.  We have the following corollary of the above
proposition.
\begin{corollary}\label{corollary:connect}
  Let $(W, a_0, a_1) \in \ddp(\point)$.  Let $U_0$ and $U_1$ be small
  open intervals in $\setR$ around $a_0$ and $a_1$, respectively,
  consisting of regular values of $f$.  Let $r: T \to \setR$ and $e:
  \coprod_T \setR^d \to W$ be as in Construction~\ref{constr:1}.  Let
  $T = T_0 \amalg T_1$ and assume $\sup U_0 < r(\tau) < \inf U_1$ for $\tau \in
  T_1$, and that
  \begin{align*}
    f^{-1}(a_1) \amalg \coprod_{T_1} \setR^d \to f^{-1}([a_0,a_1])
  \end{align*}
  is $\pi_0$\nobreakdash-surjective.  Then
  \begin{align*}
    (W^e_l, b_0, b_1) \in \ddpb(\{l\})
  \end{align*}
  for all $b_0, b_1 \in U_0 \cup U_1$ with $b_0 < b_1$, and all $l \in
  V(b_0, b_1)\cap \big(\setR^{T_0} \times \{1\}^{T_1}\big)$.
\end{corollary}
\begin{proof}
  If $b_0 \in U_0$ and $b_1 \in U_1$ then, since $U_0$ and $U_1$ are
  connected and consist of regular values of $f$,
  \begin{align}\label{eq:51}
    f^{-1}(b_1) \amalg \coprod_{T_1} \setR^d \to f^{-1}([b_0,b_1])
  \end{align}
  will also be $\pi_0$-surjective.  If $b_0, b_1 \in U_1$ or if $b_0,
  b_1 \in U_1$, then $[b_0, b_1]$ consists of regular values of $f$,
  so $f^{-1}([b_0, b_1]) \cong f^{-1}(b_1) \times [b_0, b_1]$, so the
  inclusion~(\ref{eq:51}) is $\pi_0$-surjective in this case too.
  Therefore, by Proposition~\ref{proposition:connect}(\ref{item:14})
  the element $W^{e_1} \in D_d(\setR^{T_1})$ will have
  \begin{align*}
    (W^{e_1}_{\{1\}^{T_1}}, b_0, b_1) \in \ddpb(\{1\}^{T_1}).
  \end{align*}

  It follows from Remark~\ref{remark:E-monoidal} that the construction
  of $W^e \in D_d(X \times \setR^T)$ enjoys the following naturality
  property.  If $T = T_0 \amalg T_1$, then we can restrict $e$ to
  $e_\nu: X \times \coprod_{T_\nu} \setR^d \to W$, $\nu=0,1$.  By
  construction (diagram (\ref{eq:28})), the element $W^{e_1}$ contains
  the open subset ${\proj_X^* W - \tilde e_1(X \times \setR^{T_1}
    \times \coprod_{T_1} D^d)}$ and hence $e_0$ defines an embedding
  \begin{align*}
    \proj_X^*(e_0): (X \times \setR^{T_1}) \times \coprod_{T_0}
    \setR^d \to W^{e_1}.
  \end{align*}
  The naturality property is that
  \begin{align*}
    (W^{e_1})^{\proj_X^*(e_0)} = W^e.
  \end{align*}
  Restricting to $\{1\}^{T_1} \times \setR^{T_0}$ we have
  \begin{align*}
    W^{e}_{\{1\}^{T_1}\times \setR^{T_0}} =
    (W^{e_1}_{\{1\}^{T_1}})^{\proj_X^*(e_0)}.
  \end{align*}
  The claim now follows from
  Proposition~\ref{proposition:connect}(\ref{item:15}) above.
\end{proof}

We will say that an open set $U_0
\subseteq X \times \setR$ is a \emph{tube} around $a_0$ if it contains
the graph of $a_0$, and if the intersection $U_0 \cap \{x\} \times
\setR$ is an interval consisting of regular values of $f_x: W_x \to
\setR$ for all $x \in X$.

\begin{definition}\label{definition:M-e-lambda}
  For a function $\lambda: X \times T \to [0,1]$, let $\hat \lambda :
  X \times \setR \to X \times \setR^T$ denote the adjoint $\hat
  \lambda(x,t) = (x,t\lambda(x))$.  Given $r: X \times T \to \setR$
  and $e: X \times \coprod_T \setR^d \to W$ as in
  Construction~\ref{constr:1}, let $W^{e,\lambda} \in D_d(X \times
  \setR)$ denote the pullback of $W^e$ along $\hat\lambda$.
\end{definition}

If $T = T_0 \amalg T'$ and $\lambda_{|X \times T_0} = 0$, then
$W^{e,\lambda}$ = $W^{e', \lambda'}$, where $e'$ and $\lambda'$ are
the restrictions to $T'\subseteq T$.  The following corollary follows
immediately from Corollary~\ref{corollary:connect} above.
\begin{corollary}
  \label{cor:maincor2}
  Let $(W, a_0, a_1) \in \ddp(X)$.  Let $r,e,\lambda$ be as in
  Definition~\ref{definition:M-e-lambda}.  Let $W^{e,\lambda} \in
  D_d(X \times \setR)$ be the resulting element.  Let $U_0$, $U_1$ be
  tubes around $a_0$ and $a_1$.  Assume that there is a subset $T_1
  \subseteq T$ with $\lambda_{|X \times T_1} = 1$, such that the graph
  of $r_{|X \times T_1}$ is above $U_0$ and below $U_1$, and such that
  \begin{align*}
    f_x^{-1}(a_1(x)) \amalg \coprod_{T_1} \setR^d \to
    f_x^{-1}([a_0(x), a_1(x)])
  \end{align*}
  is $\pi_0$\nobreakdash-surjective for all $x$.

  For all $b_0, b_1 : X \to \setR$ with $b_0 < b_1$ and
  $\mathrm{graph}(b_\nu) \subseteq U_0 \cup U_1$, let $\hat V(b_0,
  b_1)$ denote the intersection $X \times \{1\} \cap \hat
  \lambda^{-1}V(b_0, b_1)$.  Then the resulting element
  \begin{align*}
    (W^{e,\lambda}, b_0, b_1)_{|\hat\lambda^{-1}V(b_0, b_1)} \in
    \ddp(\hat \lambda^{-1} V(b_0, b_1))
  \end{align*}
  restricts to an element
  \begin{align*}
    (W^{e,\lambda}, b_0, b_1)_{|\hat V(b_0, b_1)} \in \ddpb(\hat V(b_0,
    b_1))
  \end{align*}
\end{corollary}
Thus, we get a \emph{concordance} from $W = W^{e,\lambda}_{|X \times
  \{0\}} \in D_d(X\times\{0\})$ to the element $W^{e,\lambda}_{|X
  \times \{1\}} \in D_d(X \times \{1\})$ and the latter element lifts over
$\hat V(b_0, b_1)$ to morphisms in $\ddpb$.

\subsection{Parametrized surgery}\label{sec:parametrized-surgery}

So far we have described how to perform surgery on $W\in D_d(X)$ along
an embedding $e: X \times \coprod_T \setR^d \to W$.  If we only have
such embeddings given \emph{locally} in $X$, then we can perform the
surgeries locally and glue them together using appropriate partitions
of unity.  More precisely we have the following construction.
\begin{construction}\label{constr:param-mul-surg}
  Let $(p,r):E \to X\times \setR$ be smooth, with $p: E \to\nobreak X$
  etale (local diffeomorphism).  Let $e: E \times \setR^d \to W$ an
  embedding over $X \times \setR$.  Let $\lambda: E \to [0,1]$ be a
  smooth map with $p|\supp{\lambda}$ proper.  Define an element
  $W^{e,\lambda} \in D_d(X \times \setR)$ in the following way.  For
  $x \in X$, the set $T_x = p^{-1}(x) \cap \supp{\lambda}$ is finite.
  Choose a connected neighborhood $U_x \subseteq X$ of $x$, and extend
  to a (unique) embedding $T_x \times U_x \to E$ over $X$, such that
  $p^{-1}(U_x) \cap \supp{\lambda}$ is contained in $T_x \times U_x$
  (this can be done because $p_{|\supp(\lambda)}$ is a closed map).

  Define $W^{e,\lambda}_{|U_x} \in D_d(U_x \times \setR)$ as the
  construction in Definition~\ref{definition:M-e-lambda} applied to
  the restriction of $e$ to $T_x\times U_x$. (If $T_x=\emptyset$ then
  $W^{e,\lambda}_{|U_x}=W_{|U_x}$.) These elements agree on overlaps, so
  by the sheaf property of $D_d$ we have defined $W^{e,\lambda} \in
  D_d(X \times \setR)$.
\end{construction}

We are now ready to prove that $\beta\ddpb \to D_d$ is a homotopy
equivalence.  It suffices to prove that any element of $D_d(X)$ is
concordant to an element which lifts to $\beta\ddpb(X)$ (plus
corresponding relative statement).

Given an element $(W, \pi, f) \in D_d(X)$, we choose (as in the proof
of Proposition~\ref{prop:5}) a locally finite open covering $X =
\cup_j E_j$ and corresponding numbers $a_j \in \setR$ such that $(W,
a_j)_{|E_j} \in \ddp(E_j)$ for all $j$.  We can assume that the $a_j$
are all distinct constants.

For each pair $j,k$ with $a_j < a_k$, let $E_{jk} = E_j \cap E_k$.
Then $\phi_{jk} = (W, a_j, a_k)_{|E_{jk}}$ is a morphism in
$\ddp(E_{jk})$.  We can assume that $E_{jk}$ is either contractible or
empty, so $(\pi,f)^{-1} (E_{jk} \times [a_j, a_k]) \cong E_{jk} \times
W_0$ for a compact manifold $W_0$ with boundary.  Consider the
inclusion
\begin{align*}
  (\pi,f)^{-1}(E_{jk} \times \{a_k\}) \to (\pi,f)^{-1}(E_{jk}\times
  [a_j, a_k]).
\end{align*}
If this is $\pi_0$-surjective, then $\phi_{jk} \in \ddpb(E_{jk})$.  If
not, we can choose a finite set $T_{jk}$ and an embedding $\tilde
e_{jk}: E_{jk} \times T_{jk} \to (\pi,f)^{-1}(E_{jk} \times (a_j,
a_k))$ over $E_{jk}$ such that
\begin{align*}
  (\pi,f)^{-1}(E_{jk} \times \{a_k\}) \amalg E_{jk} \times T_{jk} \to
  (\pi,f)^{-1}(E_{jk}\times [a_j, a_k])
\end{align*}
is $\pi_0$-surjective.  Let $r_{jk} = f \circ \tilde e_{jk}: E_{jk}
\times T_{jk} \to \setR$.  Let $E = \coprod E_{jk}\times T_{jk}$, and
let $(p,r): E \to X \times \setR$ be the resulting map.  Then the
$\tilde e_{jk}$ assemble to a map $\tilde e: E \to W$ over $X \times
\setR$.  By possibly changing the $f$\nobreakdash-level of $\tilde
e_{jk}$, we can arrange that the various $\tilde e_{jk}$ have disjoint
images so that $\tilde e$ is an embedding.  $E$ has contractible
components, so the normal bundle of $\tilde e$ can be trivialized.
Thus $\tilde e$ extends to an embedding $e: E \times \setR^d \to W$
over $X$.

Now, for each $v \in p^{-1}(x) \subseteq E$, $e$ defines an embedding
$e_v: \{v\}\times \setR^d \to W_x$, but $f_x: W_x \to \setR$ might not
be constant on the image of $e_v$. However, let $\phi:[0,\infty) \to
[0,\infty)$ be a smooth proper function with $\phi[0,1]=0$ and
$\phi'(t)>0$ for $t>1$ and $\phi(t) = t$ for $t \geq 2$. Then $f_x
\circ e_v\bigl(\phi(|u|)u\bigr)$ is constantly equal to $r(v)$ for
$u\in D^d$ and agrees with $f_x e_v(u)$ outside $2D^d$.  After
changing $f_x$ on the image of $e_v$ and then re-choosing the
embedding $e$ (precompose it with an embedding of $\setR^d$ into
$D^d$), we can assume that $e_v$ maps into $f_x^{-1}(r(v))$.  This
process works equally well in the parametrized setting, so after
modifying $f: W \to \setR$ we can assume that $e: E \times \setR^d \to
W$ is an embedding with $\pi\circ e(v,u)=p(v)$ and $f\circ
e(v,u)=r(v)$.  Choose compactly supported $\lambda_j: E_j\to [0,1]$
such that $X$ is covered by the sets $\tilde E_j = \Int
\lambda_j^{-1}(1)$, and let $\lambda_{jk} = \lambda_j\lambda_k: E_{jk}
\to [0,1]$.  These assemble to a function $\lambda: E \to \setR$ with
$p_{|\supp(\lambda)}$ proper.

Using these $p$, $r$, $e$ and $\lambda$,
Construction~\ref{constr:param-mul-surg} provides an element
$W^{e,\lambda} \in D_d(X \times \setR)$.  We claim that
$W^{e,\lambda}_1 = W^{e,\lambda}_{|X \times \{1\}}$ lifts to an
element of $\beta\ddpb(X)$.  Indeed, for $x \in \tilde E_j$, choose
$b_{xj} \in \setR$ in a tube around $a_j$ such that $(x,1) \in \hat
V(b_{xj}, b_{xj})$.  Choose a neighborhood $U_{xj}$ such that $U_{xj}
\times \{1\} \subseteq \hat V(b_{xj}, b_{xj})$.  Then
$(W^{e,\lambda}_1, b_{xj}, b_{xj})_{|U_{xj}}$ is an object of
$\ddp(U_{xj})$.  As before, refining the $U_{xj}$ to a locally finite
covering defines an element of $\beta\ddp(X)$ which in turn, by
Corollary~\ref{cor:maincor2}, is an element of $\beta\ddpb(X)$.

\begin{remark}\label{remark:Segal}
  The morphisms in Segal's cobordism category $\mathcal{S}$ are
  Riemann surfaces up to diffeomorphism.  As described in the
  introduction, an embedded oriented surface has a canonical complex
  structure (determined by being in the same conformal class as the
  euclidean metric).  Speaking loosely, this gives a functor
  $\mathcal{C}_2^+ \to \mathcal{S}$, which on morphism spaces look like
  \begin{align}\label{eq:42}
    \Emb(\Sigma, \setR^\infty)/ \Diff(\Sigma) \to
    J(\Sigma)/\Diff(\Sigma) = M(\Sigma).
  \end{align}
  Here, $J(\Sigma)$ denotes the space of complex structures on
  $\Sigma$ and $M(\Sigma)$ is the \emph{moduli space} of Riemann
  surfaces diffeomorphic to $\Sigma$.  It is a consequence of
  Teichm\"uller theory that~\eqref{eq:42} is a rational homology
  equivalence if all closed components of $\Sigma$ has genus at least
  2.  Integrally it is usually not an equivalence, due to the action
  of $\Diff(\Sigma)$ on $J(\Sigma)$ not being free in general (because
  Riemann surfaces can have non-trivial automorphisms).  From our
  point of view, it is more natural to keep track of automorphisms, by
  replacing the orbit space of the action of $\Diff(\Sigma)$ on
  $J(\Sigma)$ by the \emph{groupoid} that comes from the action.  This
  groupoid represents the \emph{moduli stack} $\mathcal{M}(\Sigma)$.
  From this point of view, $\mathcal{S}$ is a 2-category (1-morphisms
  being cobordisms with complex structure and 2-morphisms being
  isomorphisms of such), and its classifying space is homotopy
  equivalent to $B\mathcal{C}_2^+$.  We do not wish to make these
  statements precise here, or even give a precise definition of
  $\mathcal{S}$, but the main point is that the groupoid of complex
  structure on $\Sigma$ and isomorphisms of such has classifying space
  homotopy equivalent to $B\Diff(\Sigma)$ because $J(\Sigma)$ is
  contractible.

  If on the other hand we stick to the coarse moduli space
  $M(\Sigma)$, the resulting category changes, even rationally.
  Considering a 2-sphere as a cobordism from the empty 1-manifold to
  itself gives a map
  \begin{align}\label{eq:43}
    BSO(3) \simeq B\Diff^+(S^2) \to \Omega B \mathcal{C}_2^+ \simeq \setZ
    \times B\Gamma_\infty^+,
  \end{align}
  and it can be seen that the pullback of the ``Miller-Morita-Mumford
  classes'' $\kappa_i$ gives $\kappa_{2i} \mapsto 2p_i$, where $p_i
  \in H^{4i}(BSO(3))$ is the Pontrjagin class.  In particular the
  map~\eqref{eq:43} is non-trivial in rational cohomology.  If we
  replace $\mathcal{C}_2^+$ by $\mathcal{S}$, the map factors through
  the moduli space $M(S^2)$ which is a point.  Hence $\Omega
  B\mathcal{C}_2^+$ and $\Omega B\mathcal{S}$ do not have isomorphic
  rational cohomology.

  If we restrict attention to the positive boundary category, the
  difference between $\mathcal{C}_2^+$ and $\mathcal{S}$ is much less
  subject to interpretation.  Riemann surfaces with boundary cannot
  have automorphisms (which act as the identity on the boundary), and
  the map~\eqref{eq:42} is a homotopy equivalence if $\Sigma$ has no
  closed components.  Again, this can be used to prove that the
  positive boundary version of $\mathcal{S}$ has classifying space
  homotopy equivalent to $B\mathcal{C}_2^+$.
\end{remark}


\section{Harer type stability and $\mathcal{C}_2$}
\label{sec:harer-type-stability}

\cite{MR1474157} introduced a version $\mathscr{S}_b$ of the category
$\mathcal{C}_{2,\partial}^+$ to prove that $\setZ \times B\Gamma_{\infty,n}$
is homology equivalent to an infinite loop space.  This used two
properties of $\mathscr{S}_b$.  Firstly that $\mathscr{S}_b$ is
symmetric monoidal, and secondly that $\Omega B\mathscr{S}_b$ is
homology equivalent to $\setZ \times B\Gamma_{\infty,n}$.  In this
section we will prove that $\Omega B\mathcal{C}_{2,\partial}^+$ is homology
equivalent to $\setZ \times B\Gamma_{\infty,n}$, using a version of
the argument from \cite{MR1474157}.

The original stability theorem, proved by J.\ Harer in \cite{MR786348}
is about the homology of the oriented mapping class group.  In the
language used in this paper, it can be stated as follows.  Consider an
oriented surface $W_{g,n}$ of genus $g$ with $n$ boundary circles.
There are inclusions $W_{g,n} \to W_{g+1,n}$ and $W_{g,n} \to
W_{g,n-1}$ by adding the torus $W_{1,2}$ or the disk $W_{0,1}$ to one
of the boundary circles.  Let $\Diff^+(W,\partial)$ denote the group
of orientation-preserving diffeomorphisms of $W$ that restrict to the
identity near the boundary, and let
\begin{align}\label{eq:12}
  B \Diff^+(W_{g,n};\partial) &\to B \Diff^+(W_{g+1,n};\partial),\\
  \label{eq:40}
  B \Diff^+(W_{g,n};\partial) &\to B \Diff^+(W_{g,n-1};\partial),
\end{align}
be the maps of classifying spaces induced from the above inclusions.
Harer's stability theorem is that the maps in~\eqref{eq:12}
and~\eqref{eq:40} induce isomorphisms, in integral homology in a range
of dimensions that tends to infinity with $g$.  (The range is
approximately $g/2$ \cite{MR1015128}.)

In the setup of chapter~\ref{sec:tang-struct}, Harer's stability
theorem concerns the case $\theta: B \to BO(2)$, where $B = EO(2)
\times_{O(2)} (O(2)/SO(2))$.  Recently, homological stability theorems
have been proved for surfaces with tangential structure in a number of
other situations, which we now list.
\begin{itemize}
\item N.\ Wahl considered stability for non-orientable surfaces in
  \cite{Wahl06}.  Let $S_{g,n}$ denote the connected sum of $g$ copies
  of $\setR P^2$ with $n$ disks cut out, and consider the analogue
  of~\eqref{eq:12} with $\Diff^+(W_{g,n};\partial)$ replaced by
  $\Diff(S_{g,n};\partial)$.  She proves a stability range
  (approximately $g/4$) for the associated mapping class groups $\pi_0
  \Diff(S_{g,n};\partial)$ and, using the contractibility of the
  component $\Diff_1(S_{g,n};\partial)$, deduces the homological
  stability for $B \Diff(S_{g,n};\partial)$.
\item Stability for spin mapping class groups was established in
  \cite{MR1054572} and \cite{MR2068314}.  It corresponds to the
  category $\mathcal{C}_2^\theta$, with the tangential structure
  $\theta: B \Spin(2) \to BO(2)$, cf.\ \cite{Galatius06}.
\item Our final example is the stability theorem from
  \cite{CohenMadsen06}, corresponding to the tangential structure
  \begin{align*}
    \theta: EO(2) \times_{O(2)} ((O(2)/SO(2))\times Z) \to BO(2),
  \end{align*}
  where $Z$ is a \emph{simply connected} space.
\end{itemize}

With the above examples in mind, we now turn to a discussion of
abstract stability in a topological category $\mathcal{C}$.  We first
remind the reader that a square diagram of spaces
\begin{align}\label{eq:14}
\begin{aligned}
  \xymatrix{  Y \ar[r]\ar[d]^{f} & X_0 \ar[d]^g\\
    X_1 \ar[r]^p & X }
\end{aligned}
\end{align}
is \emph{homotopy cartesian} if for all $x \in X_1$ the induced map of
the vertical homotopy fibers
\begin{align}\label{eq:13}
  \hofib_{x}(f) \to \hofib_{p(x)}(g)
\end{align}
is a weak equivalence.  Similarly, the diagram~\eqref{eq:14} is
\emph{homology cartesian} if~\eqref{eq:13} is a homology equivalence,
i.e.\ induces an isomorphism in integral homology.  If the map $g$ is
a Serre fibration, then diagram~\eqref{eq:14} is homotopy cartesian if
it is cartesian.

We also remind the reader that if $\mathcal{C}$ is a category, then a
functor $F: \mathcal{C}^\mathrm{op} \to \mathrm{Sets}$ determines, and
is determined by, a category $(F \wr \mathcal{C})$ and a projection
functor $(F \wr \mathcal{C}) \to \mathcal{C}$, such that the diagram
of sets
\begin{align}\label{eq:2}
\begin{aligned}
  \xymatrix{
    {N_1(F \wr \mathcal{C})} \ar[r]^{d_i} \ar[d] &
    {N_0(F \wr \mathcal{C})} \ar[d] \\
    {N_1 \mathcal{C}} \ar[r]^{d_i} & {N_0\mathcal{C}}
  }
\end{aligned}
\end{align}
is cartesian for $i=0$ (so $d_i$ is the target map).  Explicitly, $(F
\wr \mathcal{C})$ is defined by
\begin{align*}
  N_0(F \wr \mathcal{C}) & = \{(x,c) \mid c \in N_0\mathcal{C}, x \in
  F(c)\},\\
  N_1(F \wr \mathcal{C}) & = \{(x,f) \mid f \in N_1\mathcal{C}, x \in
  F(d_0f)\}.
\end{align*}
Similarly, a functor $F$ with values in the category of spaces
determines, and is determined by, a \emph{topological} category $(F
\wr \mathcal{C})$ with a projection functor to $\mathcal{C}$ such that
the diagram~\eqref{eq:2} is a cartesian diagram of spaces for $i=0$.
If the category $\mathcal{C}$ itself is topological, then it is better
to take this as a \emph{definition}: A functor $F:
\mathcal{C}^\mathrm{op} \to \mathrm{Spaces}$ is a topological category
$(F \wr \mathcal{C})$ together with a functor $(F\wr \mathcal{C}) \to
\mathcal{C}$ such that the diagram~\eqref{eq:2} is a cartesian diagram
of spaces for $i=0$.

We return to~\eqref{eq:2} under the assumption that the right hand
vertical map is a Serre fibration.  Then the diagram is homotopy
cartesian for $i =0$.  It is homotopy cartesian also for $i=1$,
precisely if every morphism $f: x \to y$ in $\mathcal{C}$ induces a
weak equivalence $F(f): F(y) \to F(x)$.  Similarly it is homology
cartesian for $i=1$, precisely if every $f: x \to y$ induces an
isomorphism $F(f)_*: H_*(F(y)) \to H_*(F(x))$.
\begin{proposition}
  \label{prop:abstract-stability}
  Let $F: \mathcal{C}^\mathrm{op} \to \mathrm{Spaces}$ be a functor
  such that $N_0(F \wr \mathcal{C}) \to N_0\mathcal{C}$ is a Serre
  fibration.  Suppose that every $f: x\to y$ in $\mathcal{C}$ induces
  an isomorphism $F(f)_*: H_*(F(y)) \to H_*(F(x))$ and that $B(F \wr
  \mathcal{C})$ is contractible.  Then for each object $c \in
  \mathcal{C}$ there is a map
  \begin{align*}
    F(c) \to \Omega_c B\mathcal{C}
  \end{align*}
  which induces an isomorphism in integral homology.
\end{proposition}
\begin{proof}
  The assumptions imply that diagram~\eqref{eq:2} is homology
  cartesian for $i=0$ and $i=1$, and by induction every diagram of the
  form
  \begin{align*}
    \xymatrix{
      {N_k(F \wr \mathcal{C})} \ar[r]^{d_i} \ar[d] &
      {N_{k-1}(F \wr \mathcal{C})} \ar[d] \\
      {N_k \mathcal{C}} \ar[r]^{d_i} & {N_{k-1}\mathcal{C}}
    }
  \end{align*}
  is homology cartesian.  Then it follows from \cite[Proposition
  4]{MR0402733} that the diagram
  \begin{align*}
    \xymatrix{
      {N_0(F \wr \mathcal{C})} \ar[r]^{d_i} \ar[d] &
      {B(F \wr \mathcal{C})} \ar[d] \\
      {N_0 \mathcal{C}} \ar[r]^{d_i} & {B\mathcal{C}}
    }
  \end{align*}
  is homology cartesian, i.e.\ the induced map of vertical homotopy
  fibers is a homology isomorphism.  Let $c \in \mathrm{Ob}
  \mathcal{C}$.  Since $N_0(F \wr \mathcal{C}) \to N_0\mathcal{C}$ is
  assumed a Serre fibration, the homotopy fiber at $c$ of the left
  vertical map is $F(c)$.  Since $B(F \wr \mathcal{C})$ is assumed to
  be contractible, the homotopy fiber of the right vertical map at $c$
  is $\Omega_c B\mathcal{C}$.
\end{proof}

We apply this in the case where $\mathcal{C} \subseteq
\mathcal{C}_{\theta,\partial}$ is the subcategory of objects $(M,a)$
with $a < 0$, and $\theta: B \to BO(2)$ is a tangential structure for
which we have a Harer type stability theorem.  To define a functor $F:
\mathcal{C}^\mathrm{op} \to \mathrm{Spaces}$, let $S^1 \subseteq
\setR^{2-1+\infty}$ be a fixed circle, and consider the objects $b_i =
\{i\}\times S^1$ in $(\mathcal{C}_{\theta,\partial})$, $i\in \setN$.
Choose morphisms $\beta_i \subseteq [i,i+1]\times \setR^{2-1+\infty}$
from $b_i$ to $b_{i+1}$ which are connected surfaces of genus 1, and
compatible $\theta$-structures on the $b_i$ and the $\beta_i$.  We use
here that the tangent bundle of the surface $\beta_i \cong W_{1,2}$
can be trivialized.  Let $F_i: \mathcal{C}^\mathrm{op} \to
\mathrm{Spaces}$ be the functors
\begin{align*}
  F_i(c) = \mathcal{C}_{\theta,\partial}(c,b_i)
\end{align*}
and let
\begin{align*}
  F(c) = \hocolim(F_0(c) \xrightarrow{\circ \beta_0} F_1(c)
  \xrightarrow{\circ \beta_1} \cdots).
\end{align*}

As a space, $N_0(F_i \wr \mathcal{C})$ is defined by the cartesian
diagram
\begin{align*}
  \xymatrix{
    N_0(F_i \wr \mathcal{C}) \ar[r]\ar[d] &
    X_1 \ar[d]^{(d_0,d_1)}\\
    N_0 \mathcal{C} \ar[r]^-{(b_i,\mathrm{id})} &
    X_0\times
    N_0{\mathcal{C}},
  }
\end{align*}
where $X_1 = \{(W,a_0, a_1, l) \in N_1\mathcal{C}_{\theta,\partial}
\mid a_0 < 0 < a_1\}$ and $X_0 = \{(M,a,l) \in N_0
\mathcal{C}_{\theta,\partial}\mid a > 0\}$.  It follows from
\cite{MR1471480} that the right hand vertical map is a smooth Serre
fibration, so $N_0(F_i \wr \mathcal{C}) \to N_0\mathcal{C}$ and in
turn $N_0(F \wr \mathcal{C}) \to N_0 \mathcal{C}$ are Serre
fibrations, as required in Proposition~\ref{prop:abstract-stability}.
The category $(F_i \wr \mathcal{C})$ has terminal object
$\mathrm{id}_{b_i}$, so $B(F_i \wr \mathcal{C})$ is contractible.
Therefore $B(F \wr \mathcal{C}) = \hocolim_i B(F_i \wr \mathcal{C})$
is also contractible.  Finally, if $c = \{t\} \times S_n$, where $S_n
\subseteq \setR^{2-1+\infty}$ is a disjoint union of $n$ circles, then
the homotopy equivalence~\eqref{eq:9} gives
\begin{align*}
  F_i(c) \simeq \coprod_{g\geq 0} E\Diff(W_{g,n+1},\partial)
  \times_{\Diff(W_{g,n+1},\partial)} \Bun^\partial
  (TW_{g,n+1},\theta^*U_d),
\end{align*}
where $W_{g,n+1}$ is a surface of genus $g$ with $n+1$ boundary
components, and $\Diff(W_{g,n+1},\partial)$ is the topological group
of diffeomorphisms of $W_{g,n+1}$ restricting to the identity near the
boundary.

Any morphism $x \to y$ in $\mathcal{C}$ induces a map $F_i(y) \to
F_i(x)$ which corresponds to including one connected surface $W$ into
another connected surface.  After taking the limit $g \to \infty$, any
morphism $x \to y$ in $\mathcal{C}$ induces an isomorphism $H_*(F(y))
\to H_*(F(x))$ in the four cases listed above, cf.\ \cite{Galatius06},
\cite{CohenMadsen06}, \cite{Wahl06}.  In the case of ordinary
orientations we get
\begin{align*}
  F(c) \simeq \setZ \times B\Gamma_{\infty,n+1},
\end{align*}
so we get a new proof of the generalized Mumford conjecture.
\begin{theorem}[\cite{arXiv1}] There is a homology equivalence
  \begin{align*}
    \alpha: \setZ \times B\Gamma_{\infty,n} \to \Omega^\infty
    \G{2}^+.
  \end{align*}
\end{theorem}


\nocite{*}

\bibliographystyle{alpha}
\bibliography{GMTW}


\end{document}